\pdfoutput=1
\documentclass[oneside,11pt,notitlepage]{article}
\usepackage[top=2.5cm, bottom=3.5cm, left=2.0cm, right=2.0cm]{geometry}
\usepackage[toc,page]{appendix}
\renewcommand\appendixname{Appendix}
\renewcommand\appendixpagename{Appendix}
\usepackage{adjustbox}
\usepackage[utf8]{inputenc}
\usepackage{amsmath,amssymb,amsfonts,amsthm}
\usepackage{pgfplots}
\usepackage{pgfplotstable}
\usepackage{graphicx}
\usepackage{tikz}
\usepackage{makecell}
\usepackage{siunitx}
\usepackage{subcaption}
\usepackage{float}

\usepackage{mathtools}
\usepackage{enumerate} 
\usepackage{enumitem}
\usepackage{thmtools} 
\usepackage{mathtools}
\usepackage{accents}
\usepackage{titlesec}
\usetikzlibrary{arrows.meta}
\usetikzlibrary{arrows,chains,matrix,positioning,scopes,graphs,calc,patterns}
\usetikzlibrary{decorations.pathreplacing}
\usetikzlibrary{decorations.pathmorphing}
\usetikzlibrary{intersections, pgfplots.fillbetween}
\usetikzlibrary{shapes.misc}
\usetikzlibrary{tikzmark}
\tikzset{unsaturated/.style={draw=black, decorate, decoration=snake, ->}}
\tikzset{saturated/.style={draw=black, line width=1pt, rounded corners=5, ->}}
\usepgfplotslibrary{groupplots,colorbrewer, colormaps}
\usepgfplotslibrary{external}
\pgfplotsset{compat=1.18}
\usepackage{booktabs}
\usepackage{algorithm}
\usepackage{algpseudocode}
\usepackage{stmaryrd}
\usepackage{array}
\usepackage{authblk}
\usepackage{hyperref}
\usepackage{nicefrac}
\usepackage{amsmath}
\usepackage{rotating}

\DeclareMathOperator*{\argmax}{arg\,max}

\usetikzlibrary{positioning, calc, shapes, fit}
\usepgfplotslibrary{groupplots,colorbrewer, colormaps}

\newtheorem{theorem}{Theorem}

\newtheorem{proposition}{Proposition}
\newtheorem{corollary}{Corollary}
\newtheorem{definition}{Definition}

\newtheorem{example}{Example}
\newtheorem{remark}{Remark}

\newtheorem{observation}{Observation}

\usepackage{numprint}

\newcommand{\new}[1]{{\color{blue}#1}}
\newcommand\R{\mathbb{R}}
\newcommand\N{\mathbb{N}}
\newcommand\PP{\mathcal{P}}
\newcommand\TT{\mathcal{T}}
\newcommand\BB{\mathcal{B}}
\newcommand{\vertx}[1]{\mathrm{vert}(#1)}
\newcommand{\range}[1]{\llbracket #1\rrbracket}
\newcommand{\Hc}{\mathcal{H}^c}
\newcommand{\Hb}{\mathcal{H}^b}
\newcommand{\rank}[1]{\mathrm{rank}\left(#1\right)}
\newcommand{\conv}[1]{\mathrm{conv}\left(#1\right)}
\newcommand{\supp}[1]{\mathrm{supp}\left(#1\right)}
\newcommand{\bc}[1]{\mathrm{bc}(#1)}
\newcommand{\C}[1]{#1^{\mathrm{c}}}
\newcommand{\sett}[2]{\left\{#1\ \left|\ #2\right.\right\}}
\newcommand{\set}[2]{\left\{\left.#1\ \right|\ #2\right\}}

\algnewcommand\algorithmicinput{\textbf{Input:}}
\algnewcommand\INPUT{\item[\algorithmicinput]}
\algnewcommand\algorithmicoutput{\textbf{Output:}}
\algnewcommand\OUTPUT{\item[\algorithmicoutput]}
\algnewcommand\algorithmicbreak{\textbf{break}}
\algnewcommand\Break{\algorithmicbreak}

\newcommand{\todo}[1]{{\color{red}TODO {#1}}}
\pgfplotstableset{opt/.style={
		every row no 10 column no #1/.style={%
			postproc cell content/.style={
				@cell content={$\pgfmathprintnumber[fixed,precision=0]{####1}$}
			}
		}
	}
}
\pgfplotstableset{boldopt/.style={
		every row no 10 column no #1/.style={%
			postproc cell content/.style={
				@cell content={$\bf{\pgfmathprintnumber[fixed,precision=0]{####1}}$}
			}
		}
	}
}
\pgfplotstableset{boldoptstar/.style={
		every row no 10 column no #1/.style={%
			postproc cell content/.style={
				@cell content={$\bf{\pgfmathprintnumber[fixed,precision=0]{####1}}^*$}
			}
		}
	}
}

\usepackage{natbib}
\bibpunct[, ]{(}{)}{,}{a}{}{,}%
\def\bibfont{\small}%
\def\bibsep{\smallskipamount}%
\def\bibhang{24pt}%
\def\newblock{\ }%
\def\BIBand{and}%

\title{Optimization techniques for modeling with piecewise-linear functions}
\author[1,2]{P\'eter Dobrovoczki\thanks{\href{mailto:peter.dobrovoczki@sztaki.hu}{peter.dobrovoczki@sztaki.hu}}}
\author[1,2]{Tam\'as Kis\thanks{\href{mailto:tamas.kis@sztaki.hu}{tamas.kis@sztaki.hu}}}
\affil[1]{Engineering and Management Intelligence Research Laboratory, HUN-REN Institute for Computer Science and Control, Kende utca 13-17, 1111, Budapest, Hungary}
\affil[2]{Department of Operations Research, Institute of Mathematics, E\"ot\"os Lor\'and University, P\'azm\'any P\'eter s\'et\'any 1/C, 1117, Budapest, Hungary}

\date{18 February 2026}

\begin{document}
	\maketitle
	
	\begin{abstract}
		In this paper we aim to construct piecewise-linear (PWL) approximations for functions of multiple variables and to build compact mixed-integer linear programming (MILP) formulations to represent the resulting PWL function.
		On the one hand, we describe a simple heuristic to iteratively construct a triangulation with guaranteed absolute error.
		On the other hand, we extend known techniques for modeling PWLs in MILPs more efficiently than state-of-the-art methods permit.
		The crux of our method is that the MILP model is a result of solving some hard combinatorial optimization problems, for which we present heuristic algorithms.
		The effectiveness of our techniques is demonstrated by a series of computational experiments including  a short-term hydropower scheduling problem.
	\end{abstract}

	\section{Introduction}\label{sec:Intro}

	In this paper we focus on approximating a nonconvex continuous function $f\colon\ \Omega \to\R$, where $\Omega \subset \R^d$ is a bounded domain of $\R^d$ and $d \geq 1$ is an integer number, by a piecewise-linear (PWL) function $\hat{f}$, and the modeling of $\hat{f}$ in a mixed-integer linear program.
	Both topics have received considerable attention, see e.g., \cite{huchetteNonconvexPiecewiseLinear2019}, and \cite{rebennack_piecewise_2020} for recent developments.
	As an application, we reconsider the short-term hydropower scheduling problem proposed by \cite{borghettiMILPApproachShortTerm2008}.

	The input of our function approximation problem consists of a bounded rectangular domain $\Omega\subset \mathbb{R}^d$, an oracle which provides the value $f(x)$ for any $x \in \Omega$, and a parameter $\varepsilon> 0$. The output is a PWL function $\hat{f}$ such that $\max_{x \in \Omega} |f(x) - \hat{f}(x)|<\varepsilon$.
	We will use the following representation of $\hat{f}$. Let $\mathcal{T} =\{T_1,\ldots,T_m\}$ be a partitioning of $\Omega$ into $m$ simplices each of dimension $d$. Let  $S_i = \vertx{T_i}$ for each $i$,  $V = \bigcup_{i}^m S_i$ the set of all vertices, and  $\R^V_+$  the set of non-negative real vectors indexed by the vertices in $V$.
	We also define $\Delta^V = \{\lambda \in \mathbb{R}^V_+ :\ \sum_{v \in V} \lambda_v = 1 \}$, and
	$Q(S_i) := \{\lambda \in \Delta^V \ |\ \lambda_v = 0,\ \forall v \notin S_i\}$.	
	Suppose $x \in T_i$. Then there exists a unique vector $\lambda \in Q(S_i)$ such that $x = \sum_{v \in V} \lambda_v v$.
	Moreover, we define $\hat{f}(x) = \sum_{v \in V} \lambda_v f(v)$.
	While seeking a PWL function $\hat{f}$ with bounded approximation error, we try to keep the number of simplices low to limit the size of the corresponding MILP formulation.
	Our method builds on Delaunay refinement \citep{chew1989guaranteed, rupert1995delaunay, shewchuk_tetrahedral_1998} and Possion-disk sampling \citep{bridson2007fast}.

	The next problem to solve is the modeling of $\hat{f}$ in a mixed-integer linear program.
	The chosen representation of $\hat{f}$ implies that all we need is to describe the following disjunctive constraint:
	\begin{equation}
	\lambda \in \bigcup_{i=1}^m Q(S_i). \label{eq:disj_constr}
	\end{equation}
		In this paper, we extend the method based on independent branching schemes, first proposed in \cite{vielma2011modeling}. 
		 \begin{definition}[Definition 1 of \cite{vielma2011modeling}]\label{def:branching_scheme}
	The set system $\{L_\ell, R_\ell\}_{\ell=1}^D$ with $L_\ell, R_\ell \subset V$ constitutes an {\em independent branching scheme\/} of depth $D$ for constraint (\ref{eq:disj_constr}) if 
	\[
	\bigcup_{i=1}^m Q(S_i) = \bigcap_{\ell=1}^D \left(Q(L_\ell) \cup Q(R_\ell) \right).
	\]
	\end{definition}
	Independent branching schemes are also called $2-\textit{way}$ branching schemes.
	The significance of this definition is shown by the following result \cite[Theorem 4]{vielma2011modeling}. Let $\{Q(S_i)\}_{i=1}^m$ be a finite family of faces of  $\Delta^V$, and suppose there exists an independent branching scheme  $\{L_\ell, R_\ell\}_{\ell=1}^{\lceil \log_2 m\rceil}$ for (\ref{eq:disj_constr}). Then a MILP formulation for (\ref{eq:disj_constr}) is the system
	\begin{equation}
	\lambda \in \Delta^V,\quad \sum_{v\notin L_\ell} \lambda_v \leq z_\ell,\quad \sum_{v\notin R_\ell} \lambda_v \leq 1-z_\ell,\quad z_\ell \in \{0,1\}\quad \forall \ell \in \{1,\ldots,\lceil\log_2 m\rceil \},\label{eq:mip_repr}
	\end{equation}
	The existence of an independent branching scheme for a particular case of (\ref{eq:disj_constr}) is not granted. A more general concept is that of a 
	$k-\textit{way}$ independent branching scheme defined by \citet{huchette_combinatorial_2019},
	where instead of set pairs, we have set $k$-tuples in Definition~\ref{def:branching_scheme}.

	A necessary and sufficient condition for the existence of a $k-\textit{way}$ branching scheme for (\ref{eq:disj_constr}) is given by \cite[Theorem 1]{huchette_combinatorial_2019}, and is characterized in terms of the rank of the conflict hypergraph $\Hc_{\mathcal{S}}$ associated with the set system $\mathcal{S} = \{S_i\}_{i=1}^m$.
	Briefly,  the vertex set of $\Hc_{\mathcal{S}}$ is $V$, and $C\subseteq V$ is a hyperedge if no set in $\mathcal{S}$ contains $C$ as a subset, with $C$ being minimal with respect to this property.
	A major limitation is that,  although several techniques exist for constructing a 2-way branching scheme for a set system when one exists, there are no general methods for computing $k-\textit{way}$ branching schemes for $k\geq 3$.
	Motivated by this gap, we develop a new approach for representing (\ref{eq:disj_constr}) in a MILP.
	 
	The main results of this paper are the following:
	\begin{enumerate}[label=\roman*), ref=\roman*)]
		\item We have developed a general framework for finding a PWL approximation $\hat{f}$ of a function $f\colon \Omega \to\R$, where $\Omega \subset \R^d$ is a rectangular domain.
		Our method iteratively identifies a point in $\Omega$ at which the estimated absolute difference of $f$ and $\hat{f}$ is maximal, incorporates this point into the vertex set of the simplicial partition, and updates the partition accordingly.
		This process is repeated until the estimated error falls below a prescribed tolerance.
		The main challenge is to control the Lipschitz-constant of $\hat{f}$ (Section~\ref{sec:lipschitz} and Section~\ref{app:lipschitz_corr}\footnote{Some materials are presented in the Appendix consisting of Sections \ref{app:fitting_proofs}-\ref{app:blocking}.}) for which we develop a method based on repeated Delaunay refinement of a simplicial partition of $\Omega$. 
		Unlike other methods, we avoid the (explicit) solution of complicated nonlinear programs to limit the approximation error; we estimate the error by sampling $\Omega$ (Section~\ref{sec:error-est} and Section~\ref{app:sampling}) and give a guarantee that the approximation error is bounded by $\varepsilon$.
		We demonstrate the effectiveness of our method by computational experiments in Sections~\ref{sec:pwlapprox:comp} and~\ref{sec:adaptive}, and further results are presented in Section~\ref{app:fitting-exp}.
		 
		\item We propose a procedure for finding a MILP formulation for the disjunctive constraint (\ref{eq:disj_constr}) when the set system $\mathcal{S}$ corresponds to a simplicial partition of a bounded domain of $\mathbb{R}^d$ (Section~\ref{sec:mip_form} and Section~\ref{app:mip_form}). This construction immediately yields a MILP formulation for $\hat{f}$.
		The approach begins with the construction of the conflict hypergraph of $\Hc_\mathcal{S}$ associated with (\ref{eq:disj_constr}). To enable this, we prove that for any polyhedral partition of a bounded domain in $\mathbb{R}^d$, the cardinality of a minimal conflict set is at most $d+1$. This result generalizes Theorem 2 of \cite{huchette_combinatorial_2019}, and is presented in Section~\ref{sec:rank-bounding}.
		Next, we  introduce a procedure to reduce the rank of the conflict hypergraph by splitting simplices (Section~\ref{sec:rank_reduction} and Section~\ref{app:rank_reduction}). 
		After this reduction, we treat separately  the subhypegraph of $\Hc_\mathcal{S}$ induced by hyperedges of rank at least 3, and the subgraph $G^c_\mathcal{S}$ induced by rank-two edges.
		Higher-rank conflicts are resolved via a coloring-based formulation that yields additional constraints (Section~\ref{sec:high-rank-conflict} and Sections~\ref{app:coloring} and~\ref{app:higher-dim-size}), while pairwise conflicts in $G^c_{\mathcal{S}}$ are handled using a 2-way branching scheme (Section~\ref{sec:pairwise-conflicts} and Sections~\ref{app:biclique_heur} and~\ref{app:bicliqueCover:exp}). Combining the resulting formulations yields a complete MILP representation of \eqref{eq:disj_constr} (Section~\ref{sec:gib}).
	
		\item We apply our techniques to the hydropower scheduling problem of \cite{borghettiMILPApproachShortTerm2008} in Section~\ref{sec:hydropower}, and provide detailed computational results in Section~\ref{sec:experiment} and Section~\ref{app:blocking}.
	\end{enumerate}

	We conclude with a summary of the main terminology and notation.
	For a set system $\mathcal{S} = \{S_i\}_{i\in \range{m}}$ with $V = \bigcup_{i=1}^m S_i$, we say that $F \subseteq{V}$ is {\em feasible\/} if $F$ is a subset of at least one member of $\mathcal{S}$, and {\em infeasible\/} otherwise. Some $C \subseteq V$ is a {\em minimal infeasible set\/} for $\mathcal{S}$ if any proper subset of $C$ is feasible, but $C$ is not. 
	The {\em conflict hypergraph\/} for $\mathcal{S}$, denoted by $\Hc_\mathcal{S}$, has vertex set $V$, while the set of hyperedges consists of the minimal infeasible sets for $\mathcal{S}$. The {\em rank of a hyperedge\/} $E$ is its size $|E|$. The {\em rank of a hypergraph\/}  $\Hc_\mathcal{S}$ is the maximum rank of its edges.
	Let $G^c_\mathcal{S}$ denote the subgraph of 
	$\Hc_\mathcal{S}$ consisting of all the hyperedges of rank two.
	
	Given a graph $G=(V,E)$, a {\em biclique\/} $(A\cup B, E')$ in $G$ is  a complete bipartite subgraph of $G$, namely, $\emptyset \subset A,B \subset V$, $A\cap B = \emptyset$, and for all $u\in A$, and $v \in B$, $\{u,v\} \in E' \subseteq E$. A {\em biclique cover\/} of $G$ is a set of bicliques $(A_\ell\cup B_\ell, E_\ell)_{\ell=1}^K$ of $G$ such that  $E = \bigcup_{\ell=1}^K E_\ell$.
	
	A {\em polyhedral partition\/}  of a bounded domain $\Omega$ of $\mathbb{R}^d$ is a set of $d$-dimensional polytopes $\PP = \{P_1,\ldots, P_m \}$ such that $\bigcup_i P_i = \Omega$, and $P_i \cap P_j$ is a face of both $P_i$ and $P_j$ for all $P_i, P_j \in \PP$.
	The {\em vertex set of $\PP$\/} is $V_{\PP} = \bigcup_{P\in \PP} \vertx{P}$.
	When all the polyhedra in $\mathcal{P}$ are simplices,  we have a {\em simplicial partition\/}.
	
	Let $f\colon \Omega \rightarrow \R$ be a function defined on a bounded domain $\Omega \subset \R^d$, and let $\mathcal{T}$ be a simplicial partition of $\Omega$. A {\em piecewise-linear interpolation of $f$ on $\mathcal{T}$} is a PWL function $\hat{f}\colon \Omega \rightarrow \R$ such that $\hat{f}(x) = f(x)$ for all $x \in V_{\mathcal{T}}$ and $\hat{f}$ is affine on each $T \in \mathcal{T}$.
	Let $\hat{f}_T$ the be the affine function with $\hat{f}_T(x) = \hat{f}(x)$ for all $x \in T$. 
	Since $\hat{f}_T$ is an affine function, there exists $b_T \in \mathbb{R}$ such that $\hat{f}_T(x) = \nabla\hat{f}_T \cdot x + b_T$ for all $x \in T$, where $\nabla\hat{f}_T$ is the gradient of $\hat{f}_T$, a constant vector on $T$.
	 	It is  known that $\hat{f}_T$ is Lipschitz-continuous on $T$ with Lipschitz-constant $\hat{L}_T = \lVert\nabla\hat{f}_T\rVert$.
	
	We define the ball centered at $x$ with radius $r$ as $B(x,r) = \{p \in \mathbb{R}^d\colon\ \lVert x-p\rVert \leq r\}$. Let $T$ be a simplex in $\mathbb{R}^d$. We say that a finite set of points $X_T \subset T$ {\em covers $T$ with radius $r$\/} if $T \subset \bigcup_{x \in X_T} B(x,r)$.

	For convenience, we use the notation $\range{\cdot}$ for $\left[1,\cdot\right]\cap\N$ throughout the paper.

\section{Review of related literature}
\subsection{Piecewise-linear function fitting}
	For univariate functions, there are several optimization techniques for finding optimal PWL approximations, see e.g., \cite{kongDerivationContinuousPiecewise2020,rebennack_piecewise_2020}.
	For multivariate quadratic functions in $n$-dimensions, \cite{pottmann_quadratic_2000} constructed approximate PWLs, where $\hat{f}$ is specified over polyhedra that are translates of one another.
		There are several heuristics for bivariate function fitting, e.g., simulated annealing  by~\cite{schumaker_computing_1993},  a clustering-based heuristic for convex functions by~\cite{magnani_convex_2009}, and an iterative method for finding a locally optimal grid triangulation which alternates between the adjustment of the grid and the triangulation by~\cite{toriello_fitting_2012}.
		
		\citet{geisler_new_2013} propose an adaptive PWL interpolation that iteratively adds points of maximal error to the triangulation. While similar to our approach, it lacks theoretical guarantees, and \citet{burlacu_solving_2020} show that it may fail to terminate for certain functions.
		\citet{geisler_new_2013} also introduce an adaptive MILP-based refinement for gas transport optimization, refining the PWL model when constraint violations exceed a tolerance.
		A related method by \citet{burlacu_solving_2020} solves successive MILP relaxations of a MINLP, bisecting those triangles that contain an optimal point which violates a non-linear constraint, along their longest edge. This ensures convergence and yields finer triangulations near optimal points. 


\subsection{Biclique covers}
\label{sec:bicliques}
Let $G = (V, E)$ be a simple graph.
The {\em minimum biclique cover problem (MBCP)}  aims at finding a biclique cover for a graph $G$ with a minimum number of bicliques.
Let $w\colon E\to\R_{\geq0}$ be a weight function on the edges of $G$.
		The {\em maximum weight biclique problem (MWBP)} is to find a biclique $G'=(V',E')$ in $G$ such that $\sum_{e\in E'} w_e$ is maximal.
		The MWBP is \textsf{NP}-complete even if all the edge weights are equal to $1$ and the graph is bipartite~\citep {peeters_maximum_2003}.

		The recent survey of \citet{schwartzOverviewGraphCovering2022} summarizes the current state of the art regarding biclique covers; however, we list here some well-known results of the area. Let $\bc{G}$ denote the minimum size of a biclique cover of graph $G$. 
		A trivial upper bound for $\bc{G}$  is $n-1$, where $n=|V|$ (by covering the graph with $n-1$ stars).
		The best known general bound for $\bc{G}$ is $n-\lfloor\log_2{\frac{2n}{3}}\rfloor$~\citep{tuzaCoveringGraphsComplete1984}.
		In general, solving the MBCP is \textsf{NP}-hard even if $G$ is bipartite~\citep{ORLIN1977406}.
		There are mixed-integer programming formulations for the MBCP and for the MWBP as well, with an exponential number of inequalities by~\cite{cornaz_chromatic_2006}.
		A different formulation for MBCP is shown in~\cite{huchette_combinatorial_2019}.

\subsection{MILP formulations for PWL functions}

There are several approaches for modeling (non-separable) single and multivariate PWL functions in mixed-integer linear programs, see e.g.,
	the multiple choice model (MC) in~\cite{balakrishnan_composite_1989},
	the incremental model (Inc) in~\cite{wilson1998polyhedral},
	the disaggregated convex combination model (DCC) in ~\cite{meyer1976mixed},
	the convex combination model (CC) in~\cite{leePolyhedralMethodsPiecewiselinear2001},
	the disaggregated logarithmic convex combination (DLog) in~\cite{vielma2011modeling}.
	For a comprehensive study of the well-known models for MILP formulations of multivariate PWL functions, we refer the reader to~\cite{vielmaMixedIntegerModelsNonseparable2010}.
	\cite{huchette_combinatorial_2019} generalize some of the results of \cite{vielma2011modeling}, and they also prove that for polyhedral partitions in $\mathbb{R}^2$, the maximum size of a minimal conflict set is 3.
	The 6-stencil formulation \citep{huchetteNonconvexPiecewiseLinear2019} offers a systematic way of constructing a pairwise-independent branching scheme for grid triangulations, while \cite{lyu2024modeling} provide a heuristic for finding a biclique cover when the underlying combinatorial disjunction admits a junction tree.

	\section{The procedure for computing a MILP-approximation of a nonlinear function}
	\label{sec:method_overview}

	In this section, we present an overview of our method, which takes as input a nonlinear function $f$, and produces as output a  MILP formulation of a PWL approximation $\hat{f}$ of $f$. We assume that $f\ :\ \Omega \longrightarrow \R$, where $\Omega \subset\mathbb{R}^d$ is a rectangular domain for some integer $d\geq 1$. 
	
	Our method has two main components. The first  constructs a PWL approximation $\hat{f}$ of $f$, and the second  derives a MILP formulation representing $\hat{f}$.
	The combined workflow of these two components is illustrated in the left and right panels of Figure~\ref{fig:flowchart}, respectively.

		Our procedure for constructing $\hat{f}$ begins with a simplicial partition of $\Omega$, whose vertices coincide with the corner points of $\Omega$.
		This initial partition defines a PWL interpolant $\hat{f}$ of $f$, which is iteratively refined together with $\hat{f}$ -- until a sufficiently accurate approximation of $f$ is obtained.
		In each iteration, first a Lipschitz-correction step is applied to enforce a bound on the Lipschitz-constant of $\hat{f}$ (Section~\ref{sec:lipschitz}).
		Subsequently, the maximal approximation error is estimated  by sampling a suitably sized set of points from the simplices defining $\hat{f}$ (Section~\ref{sec:error-est}).
		If the estimated error exceeds the prescribed bound  $\epsilon$,   the simplicial partition is extended by incorporating the point of maximal error, and the method continues with Lipschitz-correction. Otherwise, the construction of $\hat{f}$ terminates. See Algorithm~\ref{alg:fitting} for a formal description.

After computing $\hat{f}$, the method proceeds with the construction of the conflict hypergraph $\Hc_{\mathcal{S}}$ (Section~\ref{sec:rank-bounding}).
It then attempts to reduce the rank of $\Hc_\mathcal{S}$ by modifying the simplicial partition of $\Omega$ (Section~\ref{sec:rank_reduction}).
If, after the reduction step, the rank of $\Hc_{\mathcal{S}}$ remains at least three,  the blocking hypergraph associated with $\mathcal{S}$ is constructed, and a vertex coloring  is computed such that no hyperedge is monochromatic (Section~\ref{sec:high-rank-conflict}).
For the subhypergraph of $\Hc_{\mathcal{S}}$  spanned by rank-two edges, our method finds a biclique cover (Section~\ref{sec:pairwise-conflicts}).
Using the coloring and the biclique cover, it finally constructs the MILP formulation for (\ref{eq:disj_constr}) in Section~\ref{sec:gib}.

	\begin{figure}
	\tikzstyle{startstop} = [rectangle, rounded corners, minimum width=3cm, minimum height=1cm,text centered, draw=black]
	\tikzstyle{decision} = [diamond, draw, text width=6.5em, text badly centered, inner sep=0pt]
	\tikzstyle{process} = [rectangle, minimum width=3cm, minimum height=1cm, text centered, draw=black]
	\tikzstyle{arrow}  = [thick,->,>=stealth]
		\resizebox{1\textwidth}{!}{
			\tikzexternaldisable
			\begin{tikzpicture}[node distance = 1.5cm]					
				\node (start) [startstop, align=center] {Input:\\ $f,\varepsilon, \Omega$};
				\node (init) [process, below of=start, yshift=-.5cm, align=center] {Initialize ground points $V$ \\  and simplicial partitioning $\TT$};
				\node (refine) [process, below of=init, yshift=-.5cm, align=center] {Lipschitz-correction \\ (Section~\ref{sec:lipschitz})};
				\node (estimate) [process, below of=refine, yshift=-.5cm, align=center] {Error estimation \\ (Section~\ref{sec:error-est})};
				\node (tol) [decision, below of=estimate, yshift=-1.5cm] {Is the error below $\varepsilon$?};
				\node (triangulation) [process, right of=estimate, xshift=4cm, align=center] {Determine set system $\mathcal{S}$ \\ from $\TT$};
				\node (newpt) [process, left of=estimate, xshift = -4cm, align=center] {Extend $V$ with a new point of\\ maximum error and update $\mathcal{T}$};
				\node (constructConflict) [process, right of=start, xshift=8cm, align=center] {Construct conflict hypergraph $\Hc_{\mathcal{S}}$ \\ (Section~\ref{sec:rank-bounding}) };
				\node (subdivide) [process, below of=constructConflict, yshift=-1cm, align=center] {Reducing the rank of $\Hc_{\mathcal{S}}$\\ (Section~\ref{sec:rank_reduction})};
				\node (isRankGreater) [decision, right of=subdivide, xshift=4.5cm] {Is $\mathrm{rank}(\Hc_{\mathcal{S}}) > 2$?};
				\node (constructBlocking) [process, right of=isRankGreater, xshift=4.5cm, align=center] {Construct blocking hypergraph $\Hb_{\mathcal{S}}$ \\(Section~\ref{sec:high-rank-conflict})};
				\node (color) [process, below of=constructBlocking, yshift=-2cm, align=center] {Color blocking hypergraph\\ (Section~\ref{app:coloring})};
				\node (bicliqueCover) [process, below of=isRankGreater, yshift=-2cm, align=center] {Find biclique cover of $G^c_{\mathcal{S}}$ \\(Section~\ref{sec:pairwise-conflicts})};
				\node (constructMIP) [process, below of = bicliqueCover, yshift=-.5cm, align=center] {Construct MILP formulation\\ (Section~\ref{sec:gib})};
				\node (stop) [startstop, right of=constructMIP, xshift=4.5cm, align=center] {Output: MILP representation of \\ PWL approximation of $f$};
				\draw [arrow] (start) -- (init);
				\draw [arrow] (init) -- (refine);
				\draw [arrow] (refine) -- (estimate);
				\draw [arrow] (estimate) -- (tol);
				\draw [arrow] (tol) -| node[near start, above] {No} (newpt);
				\draw [arrow] (newpt) |- (refine);
				\draw [arrow] (tol) -| node[near start, above] {Yes} ++(5.5cm, 0) -- (triangulation.south);
				\draw [arrow] (triangulation) |- 	(constructConflict.west);
				\draw [arrow] (constructConflict) -- (subdivide);
				\draw [arrow] (subdivide) -- (isRankGreater);
				\draw [arrow] (isRankGreater) -- node[near start, left] {No} (bicliqueCover);
				\draw [arrow] (isRankGreater) --node[above] {Yes} (constructBlocking);
				\draw [arrow] (constructBlocking) -- node[right] {No} (color);
				\draw [arrow] (color) -- (bicliqueCover);
				\draw [arrow] (bicliqueCover) -- (constructMIP);
				\draw [arrow] (constructMIP) -- (stop);
				\node[
					draw,
					rounded corners=6pt,
					fit=(start) (tol) (newpt) (init),
					label={[anchor=north west]north west:{Section~\ref{sec:pwl_fit}}}
				] {};	
				\node[minimum height=0pt, minimum width=0pt] (tolshift) at ($(tol.south)+(4cm,.15cm)$) {};			
				\node[
					draw,
					rounded corners=6pt,
					fit=(triangulation) (stop) (constructBlocking) (constructConflict) (tolshift),
					label={[anchor=north east]north east:{Sections~\ref{sec:mip_form}}}
				] {};
			\end{tikzpicture}
			\tikzexternalenable
		}
	\caption{Flowchart of the PWL function fitting and the MILP formulation procedures.\label{fig:flowchart}}
\end{figure}

		\section{Piecewise-linear function fitting}
		\label{sec:pwl_fit}
		In this section, a new iterative heuristic approach is proposed to  approximate a Lipschitz continuous nonlinear function $f\colon\R^d\to\R$ by a PWL function $\hat{f}$ on a bounded rectangular domain $\Omega \subset \R^d$.
		To this end, we build iteratively a simplicial partitioning $\mathcal{T}$ of $\Omega$ which uniquely determines a PWL interpolant $\hat{f}$ of $f$ on $\Omega$. The objective is to ensure that the maximal absolute difference $\max_{x\in\Omega}|f(x)-\hat{f}(x)|$ of $f$ and $\hat{f}$ on $\Omega$ is less than a given $\varepsilon>0$.

		Our method is described formally in Algorithm~\ref{alg:fitting}. The main steps are already explained in Section~\ref{sec:method_overview}, and they are further elaborated in Sections~\ref{sec:lipschitz} and \ref{sec:error-est}. Our computational results are summarized in Section~\ref{sec:pwlapprox:comp}.
		For detailed proofs of all statements of this section refer to Section~\ref{app:fitting_proofs}.
		\begin{algorithm}
			\small
			\caption{Randomized piecewise-linear function fitting}\label{alg:fitting}
			\begin{algorithmic}[1]
				\INPUT Lipschitz-continuous function $f$ with Lipschitz-constant $L$, $\varepsilon>0$ error tolerance, polyhedral domain $\Omega\subset\R^d$.
				\OUTPUT PWL function $\hat{f}$ s.t. $\max_{x \in \Omega}|f(x)-\hat{f}(x)|\leq \varepsilon$.
				\State $V=\mathrm{ext}(\Omega),\, \mathcal{T} = \mathrm{Delaunay}(V)$ \Comment{Initialization}
				\While{true}
				\State $\mathcal{T}\gets\mathrm{Lipschitz-correction}(f,L,\mathcal{T})$ \Comment{Lipschitz-correction}
				\State $X_T = \mathrm{Sampling}(T,\nicefrac{\varepsilon}{2(L + \hat{L}_T})\quad \forall T\in\mathcal{T}$\Comment{Sampling}
				\State $p_{\max} \gets \argmax_{x\in X_T,\,T\in\mathcal{T}}|f(x) - \hat{f}(x)|$,\quad $\hat{\varepsilon}_{\max} = |f(p_{\max}) - \hat{f}(p_{\max})|$ \Comment{Error estimation}
				\State {\bf if} $\hat{\varepsilon}_{\max} \leq \varepsilon/2$ {\bf then quit the loop}\Comment{Termination condition}
				\State $V\gets V\cup\left\{p_{\max}\right\},\, \mathcal{T}\gets \mathrm{Delaunay}(V)$\Comment{Error-improvement}
				\EndWhile
				\State\Return $\hat{f}$
			\end{algorithmic}
		\end{algorithm}
		
		
		\subsection{Lipschitz-correction}\label{sec:lipschitz}
				Given a Lipschitz-continuous function $f\colon \Omega\longrightarrow\R$ with Lipschitz-constant $L$. A \emph{Lipschitz-correction algorithm for $f$}  transforms any PWL interpolant $\hat{f}$ of $f$ into one with Lipschitz-constant $\hat{L} \leq c \cdot L$, where $c$ is a fixed constant that may depend of $f$ only.

		For any  $T \in \mathcal{T}$, let $\ell^{\max}_T$ denote the length of the longest edge of simplex $T$ and $\delta^{\min}_T$  the shortest distance (in Euclidean norm) between two disjoint faces of $T$.
		First we establish a connection between the Lipschitz-constant of $f$ and that of $\hat{f}$. Recall the definition of $\hat{f}_T$ from Section~\ref{sec:Intro}.

			\begin{theorem}\label{thm:lipschitz-bound}
				Let $f\colon\Omega\to\R$ be a Lipschitz-continuous function with Lipschitz-constant $L$, $\mathcal{T}$  a simplicial partitioning of $\Omega\subset\R^d$ and  $\hat{f}$ the corresponding PWL interpolation of $f$. Then
				\begin{equation}
					\lVert\nabla\hat{f}_T\rVert\leq L\cdot\nicefrac{\ell^{\max}_T}{\delta^{\min}_T}\quad \forall T\in\mathcal{T}. \label{eq:lip-bound}
				\end{equation}
			\end{theorem}

			\begin{corollary}
				If $f\colon\Omega\to\R$  Lipschitz-continuous with Lipschitz-constant $L$, then the PWL interpolation $\hat{f}$ of $f$ is also Lipschitz-continuous with a Lipschitz-constant $\hat{L}\leq L\cdot\max_{T\in \mathcal{T}} \nicefrac{\ell^{\max}_T}{\delta^{\min}_T}$.
			\end{corollary}

			For $d=2$, the right hand side of bound~(\ref{eq:lip-bound}) evaluates to $L\cdot\nicefrac{\ell^{\max}_T}{h_T^{\min}}$ where $h_T^{\min}$ denotes the minimal height of triangle $T$.
			In fact, there is a slightly stronger bound if $\Omega \subset \mathbb{R}^2$.

		\begin{theorem}\label{thm:max-lipschitz}
			Let $\mathcal{T}$ be a triangulation of $\Omega\subset\R^2$ and $f\colon\Omega\to\R$ a Lipschitz-continuous function with Lipschitz-constant $L$. Let $\hat{f}_T$ be the affine interpolant of $f$ over $T\in\mathcal{T}$. Then
			\begin{equation}
				\lVert\nabla\hat{f}_T\rVert\leq \frac{L}{\sin\left(\alpha_T^{\min}\right)},
			\end{equation}
			where $\alpha_T^{\min}$ denotes the smallest angle of triangle $T$.			
		\end{theorem}
		\begin{corollary}
			For $d=2$, the affine interpolation $\hat{f}$ of $f$ is also Lipschitz-continuous with Lipschitz-constant $\hat{L}\leq\nicefrac{L}{\sin\left(\alpha^{\min}_{LB}\right)}$, where $\alpha^{\min}_{LB}$ is a lower bound on the  smallest angle of any triangle in the underlying triangulation of $\hat{f}$.
		\end{corollary}
				
		In $\mathbb{R}^2$ there are various methods to increase the minimum angle of a triangulation based on Delaunay refinement \citep{rupert1995delaunay,chew1989guaranteed}.
		We use the algorithm of \cite{rupert1995delaunay} which ensures that the minimum angle in the triangulation is greater than $\alpha_{\mathrm{LB}}$, if $\alpha_{\mathrm{LB}}$ is chosen to be less than $20^\circ$, refer to Algorithm~\ref{app:ruppert_alg} of Section~\ref{app:pseudocodes_subroutines}.
			By Theorem~\ref{thm:max-lipschitz}, bounding the minimum angle of the triangles from below imposes an upper bound of $c\cdot L$ on $\hat{L}$ with $c=\nicefrac{1}{\sin\left(\alpha_{\mathrm{LB}}\right)}$.

		In $\mathbb{R}^d$ for $d \geq 3$, Delaunay refinement or weighted Delaunay-triangulations are shown to have a potential in eliminating badly shaped simplices in practice. However, they lack theoretical guarantees, refer to \citep{shewchuk_tetrahedral_1998,cheng_sliver_2000, si2008three}.
		Also, one may exploit the properties of the approximated function (e.g. convexity, monotonicity, etc.) when designing a Lipschitz-correction algorithm for a specific use case.
		Note that it is not necessary for a Lipschitz-correction algorithm to add new vertices and simplices to the partitioning, it can also work by perturbing the geometric positions of the points in $V$.

	\subsection{Error estimation with sampling}\label{sec:error-est}
	In this section we describe our method of error estimation based on sampling. The main idea is that in each simplex $T \in \mathcal{T}$ we generate a finite set of points $X_T$ covering $T$ with radius $r_T = \nicefrac{\varepsilon}{2(L +\hat{L}_T)}$ (for definitions refer to Section~\ref{sec:Intro}). We will prove that
	the maximum absolute difference $\varepsilon_T = \max_{x\in T}|f(x)-\hat{f}(x)|$ of $f$ and $\hat{f}$ on any simplex  $T \in\mathcal{T}$ can be bounded by a linear function of $\hat{\varepsilon}_T = \max_{x\in X_T}|f(x)-\hat{f}(x)|$ and the Lipschitz-constants $L$ of $f$, and $\hat{L}_T$ of $\hat{f}_T$. After elaborating on choosing $X_T$ by sampling $T$, we prove our main result about the convergence of Algorithm~\ref{alg:fitting}.
  	
			\begin{proposition}\label{prop:sample-error}
			If $X_T \subset T$ covers $T$ with radius $r_T$, then $\varepsilon_T \leq \hat{\varepsilon}_T+ (L+\hat{L}_T)r_T$.
		\end{proposition}

	Based on Proposition~\ref{prop:sample-error}, we generate a finite point set $X_T \subset T$ which covers $T$ with radius $r_T = \nicefrac{\varepsilon}{2(L +\hat{L}_T)}$.
	\begin{remark}
	Generating samples with a radius smaller than $\nicefrac{\varepsilon}{2(L +\hat{L}_T)}$ leads to more accurate error estimates and, consequently, permits a looser termination criterion; however, this comes at the cost of requiring a larger number of sample points.
	By selecting $r_T=\nicefrac{\theta\varepsilon}{(L+\hat{L}_T)}$ for some $0< \theta< 0.5$, the termination condition can be  relaxed to $\hat{\varepsilon}_{\max}\leq (1-\theta)\varepsilon$.
	\end{remark}

	To find $X_T$, we adopt the maximal Poisson-disk sampling procedure of \cite{ebeida_simple_2012} (depicted in Algorithm~\ref{alg:mps} in Section~\ref{app:pseudocodes_subroutines}), for its radius control and guaranteed full coverage of the sampled domain, which works in any dimensions. 
	For efficiency, only newly formed simplices are sampled, and the point with maximum error is retained.

	The following is the main results of this section.
	
	\begin{theorem}\label{thm:algorithm}
	Algorithm~\ref{alg:fitting} computes a PWL function $\hat{f}$ with $\max_{p\in\Omega}|f(p)-\hat{f}(p)|\leq\varepsilon$ within a finite number of iterations provided the following two conditions hold:
		\begin{enumerate}
		\item there exists a Lipschitz-correction algorithm for $f$, and
		\item for each simplex $T\in \mathcal{T}$, $X_T$ is always chosen such that $X_T$ covers $T$ with radius  $r_T = \nicefrac{\varepsilon}{2(L+\hat{L}_T)}$.
		\end{enumerate}
	\end{theorem}
	
	If the conditions of this theorem do not hold, then we can still stop the algorithm at any iteration, but we have no guarantee for bounded error.

			\subsection{Computational evaluation}\label{sec:pwlapprox:comp}
			We evaluated Algorithm~\ref{alg:fitting} on various test functions. 
			Figure~\ref{fig:ripple} depicts a representative result. 
			The left panel plots the maximum error per iteration, which -- despite temporary increases due to re-triangulation -- converges to zero. 
			The right panel depicts the adaptive mesh overlaid on the target function, with finer refinement in regions of higher complexity and coarser resolution elsewhere. 
			Additional examples are presented in Section~\ref{app:fitting-exp}.

		\begin{figure}[h]

			\caption{PWL interpolation of $\sin(50\sqrt{x^2+(y-\nicefrac{1}{2})^2)}\cdot{\rm e}^{-10(x^2+(y-\nicefrac{1}{2})^2)}$ on the $\left[0,1\right]^2$ domain with absolute error $\leq 0.05$.\label{fig:ripple}}
		\end{figure}
		
	\section{A MILP-formulation for}
	\label{sec:mip_form}
	
	Suppose the PWL approximation $\hat{f}$ of $f$ is defined by the simplicial partitioning $\TT = (T_i)_{i=1}^m$ of $\Omega$, where each $T_i$ is a simplex in $\mathbb{R}^d$ as described in Section~\ref{sec:Intro}. To obtain a MILP-formulation for $\hat{f}$, we have to find a MILP-formulation for constraint (\ref{sec:Intro}).
	
	The procedure for constructing a MILP-formulation for $\hat{f}$ is explained in Section~\ref{sec:method_overview}, and its main steps are elaborated in Sections~\ref{sec:rank-bounding}-\ref{sec:gib}.

	\subsection{Constructing the conflict hypergraph $\Hc_\mathcal{S}$}	\label{sec:rank-bounding}
	To construct the blocking hypergraph for a simplicial partition $\mathcal{S}$, we need to know the rank of $\Hc_\mathcal{S}$, i.e., the maximum size of a hyperedge. 
	The following result generalizes Theorem 2 of \cite{huchette_combinatorial_2019}, and provides a general bound for polyhedral partitions in $\mathbb{R}^d$.

	\begin{theorem}\label{thm:conflict-hg}
		For  a polyhedral partition $\PP$ in $\mathbb{R}^d$,  the rank of the conflict hypergraph is at most $d+1$.
	\end{theorem}
	\begin{proof}{Proof}
		Suppose $\PP = (P_i)_{i=1}^m$ is a polyhedral partition of a bounded domain of  $\mathbb{R}^d$.
		Let $S_i = \vertx{P_i}$ ($i\in\range{m}$), and $V = \bigcup_{i=1}^{m} S_i$.
		Let $\mathcal{S} = \{ S_i\ |\ i\in \range{m}\}$.
		We prove that the maximum size of a minimal infeasible set for $\mathcal{S}$ is at most $d+1$.
		
		Suppose there exists a subset of $n \geq d+2$ points $T = (v_j)_{j=1}^{n} \subseteq V$, such that each $(n-1)$-element subset $T_j = T \setminus \{v_j\}$ of $T$ is feasible, i.e., $T_j\subseteq S_{i_j}$ for some $i_j\in \range{m}$, but $T$ is infeasible. 
		Clearly, $i_j \neq i_{j'}$ for any pair of distinct $j\neq j' \in \range{n}$, otherwise, if $i_j = i_{j'}$ for some $j\neq j'$, then $T = T_j \cup T_{j'} \subseteq S_{i_j}$ and it follows that $T$ is a feasible subset of $V$, a contradiction.
		
		Let $P'_j := \conv{T_j}$ for $j \in \range{n}$. Then $P'_j \subseteq P_{i_j}$, since $T_j \subseteq S_{i_j}$ by the definition of the index $i_j$, and $S_{i_j} \subset P_{i_j}$, since $\vertx{P_{i_j}} = S_{i_j}$ by definition. Any subset of $d+1$ polytopes from $\{P'_1,\ldots,P'_n\}$ has a point in common, since $n \geq d+2$, $|T| = n$, and each $T_j$ misses precisely one element of $T$.  Then, Helly's theorem for convex polytopes \citep{danzer1963helly} implies that  there exists some $x \in \bigcap_{j=1}^{n} P'_j$. Then $x$ is a common point of all the $P_{i_j}$.
		Since $\PP$ is a polyhedral partition, $x$ is on a common face of the polytopes $P_{i_j}$, $j \in \range{n}$.
		
		Since $x \in \bigcap_{j=1}^{n} P'_j$, $x$ can be expressed as a convex combination of the points in $T$.
		Hence, there exist non-negative coefficients $\lambda_v$ $(v \in T)$ such that $x = \sum_{v \in T} \lambda_v v$ and $\sum_{v\in T} \lambda_v = 1$. By Caratheodory's theorem, $x$ is  a convex combination of at most $d+1$ points among the points of $T$. Let $I \subset T$ be the support of $x$, i.e., $I = \{v \in T\ |\ \lambda_v > 0\}$.
		Clearly, $I$ is not empty, and has at most $d+1$ elements.
		
		We have $I\subseteq  T_j$ for some $j\in \range{n}$, since the $T_j$ constitute all the $(n-1)$-element subsets of $T$.
		Since $x$ is on a common face $F$ of the polytopes $P_{i_{j'}}$ ($j'\in \range{n}$), 
		and $I\subseteq T_j$ for some $j$, it follows that all the points of $I$ are on $F$, since $T_j$ is a subset of  $\vertx{P_{i_j}}$.
		Hence, at least one of the points $v_j \in T$ belongs to all of the polytopes $P_{i_{j'}}$ ($j'\in \range{n}$).
		Consequently,  $v_j \in S_{i_j}$ and since $T_j \subset S_{i_j}$, we have  $T = \{v_j\} \cup T_j \subseteq S_{i_j}$, a contradiction.
	\end{proof}
	
	Now, we give an example showing that the condition of the theorem is also necessary.
	\begin{example}
		In  $\mathbb{R}^d$ consider a unit sphere $S^{d-1}$ and $m \geq d+2$ points $(p_i)_{i=1}^m$ in general position (i.e., no $d+1$ points are on a common hyperplane) on $S^{d-1}$.
		The convex hull of any subset of $m-1$ points is a convex polytope $P_i$, for $i\in\range{m}$. 
		Notice that each $P_i$ is a $d$-dimensional polytope by our assumption.
		Let $\PP$ be the family of these $m$ convex polytopes.
		Then, any subset of $m-1$ points is feasible, but all the $m$ points are not, since none of the polytopes in $\PP$ contains all of them.
		Observe that  $\PP$ is not a polyhedral partition, since any pair of the polytopes in $\PP$ has a common interior point. 
	\end{example}
	We use Theorem~\ref{thm:conflict-hg} for constructing the conflict hypergraph by enumerating only those subsets of $V$ with cardinality at most $d+1$, and keeping only those that constitute a minimal conflict set.
	\subsection{Reducing the rank of $\Hc_{\mathcal{S}}$ by simplex splitting} \label{sec:rank_reduction}
	In this section we describe a procedure to reduce the rank of $\Hc_{\mathcal{S}}$ by splitting some simplices in the simplicial partitioning $\mathcal{T}$ of $\Omega$.
	The general idea is to introduce new conflicts of rank 2 between two vertices in a high rank conflict.
	This is done by splitting an edge of a simplex by a new vertex, and consequently, splitting all simplices incident to that edge into two simplices.
	
	Consider a simplicial partitioning $\mathcal{T}$  of $\Omega\subset\R^d$, and the associated set system $\mathcal{S}$ as defined in Section~\ref{sec:Intro}.
	Pick a hyperedge $C$ of $\Hc_{\mathcal{S}}$ of rank at least 3, and let $u,v$ be distinct vertices of $C$.
	Since $C$ is a minimal conflict set, $(u,v)$ is an edge of at least one simplex in $\mathcal{T}$. Let $\mathcal{S}_{uv} = \{S \in \mathcal{S}\colon\ u,v \in S\}$, and $\mathcal{T}_{uv}$ consist of those simplices containing both  $u$ and $v$.
	Clearly, $S \in \mathcal{S}_{uv}$ if and only if $\conv{S} \in \mathcal{T}_{uv}$.
	We can eliminate conflict set $C$ by the following transformation.
	We split edge $(u,v)$ by a new vertex $w \notin V$, which means that we replace each $S \in \mathcal{S}_{uv}$ with two sets
	$S_u = S\setminus \{u\} \cup \{w\}$ and $S_v = S\setminus \{v\} \cup \{w\}$.
	Accordingly, we replace $T = \conv{S}$ by the simplices $T_u = \conv{S_u}$ and $T_v = \conv{S_v}$.
	Finally, we add $w$ to $V$. After the split, $\{u,v\}$ is a conflict set, since no set of $\mathcal{S}$ (no simplex of $\mathcal{T}$) contains both $u$ and $v$. 

		\begin{observation}\label{prop:splitting}
			Splitting edge $(u,v)$  eliminates each hyperedge $C$ of $\Hc_{\mathcal{S}}$ of rank at least 3 that contains both of the vertices $u$ and $v$, since $C$  no-longer represents a minimal conflict set. 
		\end{observation}

	The following properties of the splitting operation establishes the basis for a heuristic for the reduction of the rank of the conflict hypergraph. Detailed proofs can be found in Section~\ref{app:rank_reduction}.

	\begin{proposition}\label{prop:splitting2}
		Suppose $C$ is a new minimal conflict set after splitting edge $(u,v)$ with vertex $w$, and $|C|\geq 3$. Then $w\in C$ and one of the following three cases holds before the split:
		\begin{enumerate}[label=\roman*), ref=\roman*)]
			\item $C\setminus\{w\} \cup \{u\}$ is a minimal conflict set,\label{split:case1}
			\item $C\setminus\{w\} \cup \{v\}$ is a minimal conflict set,\label{split:case2}
			\item $C\setminus\{w\} \cup \{u,v\}$ is a minimal conflict set.\label{split:case3}
		\end{enumerate} 
	\end{proposition}

	\begin{proposition}\label{prop:splitting3}
		Suppose $C$ is a new minimal conflict set after splitting edge $(u,v)$ with new vertex $w$, and $|C|\geq 3$. Then any $w'\in C\setminus\{w\}$ is in a common simplex with  $u$, and also with $v$.
	\end{proposition}

	We close this section with a heuristic to reduce the rank of $\Hc_{\mathcal{S}}$, namely, Algorithm~\ref{alg:rank-reduction} iteratively tries to reduce the rank of $\Hc_{\mathcal{S}}$ by splitting simplices in $\mathcal{T}$ by their edges. 
	That is, let $k$ be the rank of the actual $\Hc_{\mathcal{S}}$, and $\mathcal{E}_k = \left\{(u,v)\colon\exists C\in\Hc_\mathcal{S} \text{ such that } (u,v)\subset C \text{ and } |C|\geq k \right\}$.
	For each edge $(u,v)\in \mathcal{E}_k$, let $r_{uv}^k$ and $c_{uv}^k$ be the number of conflict sets of rank $k$ eliminated and created by splitting $(u,v)$, respectively.
	It selects an edge $(u^*,v^*) \in \mathcal{E}_k$ for splitting which maximizes the difference $\delta = r_{uv}^k - c_{uv}^k$.
	If $\delta \leq 0$, the algorithm stops.
	Otherwise, it splits edge $(u^*,v^*)$ with its midpoint in  $\mathbb{R}^d$, and updates the simplicial partitioning $\mathcal{T}$ and $\Hc_{\mathcal{S}}$.
	
		\begin{algorithm}
			\small
			\caption{Rank reduction heuristic for conflict hypergraph $\Hc_{\mathcal{S}}$}\label{alg:rank-reduction}
			\begin{algorithmic}[1]
				\INPUT Simplicial partitioning $\mathcal{T}$ and the associated conflict hypergraph $\Hc_{\mathcal{S}}$ with rank $k\geq 3$
				\OUTPUT Updated simplicial partitioning $\mathcal{T}$ and hypergraph $\Hc_{\mathcal{S}}$ with rank $k'\leq k$
				\While{$k\geq 3$}
				\State $\mathcal{E}_k\gets \left\{(u,v)\colon\exists C\in\Hc_\mathcal{S} \text{ such that } (u,v)\subset C \text{ and } |C|\geq k \right\}$
				\State Compute $r^k_{uv}$ and $c^k_{uv}$ for each edge $(u,v)\in\mathcal{E}_k$
				\State $(u^*,v^*)\gets\arg\max_{(u,v)\in\mathcal{E}_k}r^k_{uv}-c^k_{uv},\ \delta\gets r^k_{u^*v^*}-c^k_{u^*v^*}$ 
				\State {\bf if} $\delta \leq 0$ {\bf then quit }  the {\bf while} loop
				\State Split edge $(u^*, v^*)$ with its midpoint and update $\mathcal{T}$, $\Hc_{\mathcal{S}}$ and its rank $k$
			\EndWhile
			\State\Return $\mathcal{T}$ and $\Hc_{\mathcal{S}}$
			\end{algorithmic}
		\end{algorithm}
		
	The value of $r^k_{uv}$ can be computed by counting those minimal conflict sets of rank $k$ that contain $u$ and $v$.
	By Observation~\ref{prop:splitting}, those sets are no longer minimal conflict sets after splitting edge $(u,v)$.
	Computing $c^k_{uv}$ involves counting those minimal conflict sets of rank $k$ that satisfy one of \ref{split:case1} or \ref{split:case2} of Proposition~\ref{prop:splitting2} and also satisfy Proposition~\ref{prop:splitting3}.
	
	A crucial consideration in the application of the simplex splitting algorithm for a  PWL function $\hat{f}$  is the choice of the value of the  function at the new points. 
	A natural choice is to preserve the function value  before splitting.

\subsection{Resolving higher-rank conflicts by coloring a hypergraph}\label{sec:high-rank-conflict}
	In this section, we address the conflicts corresponding to hyperedges of  $\Hc_{\mathcal{S}}$ with rank at least 3.
		We say that a subset $B$ of sets from $\mathcal{S}$ \emph{spans a minimal infeasible set $C$ of  $\Hc_\mathcal{S}$\/} if $C \subseteq \bigcup_{S \in B} S$.
		We say that $B \subseteq \mathcal{S}$ is \emph{blocking}
		if it spans a minimal infeasible set of size at least 3, and it is \emph{minimal blocking} if any proper subset of it is non-blocking. The detailed proofs of the statements of this Section can be found in Section~\ref{app:high-rank-conflict}.
		
		\begin{definition}\label{def:blocking} The blocking hypergraph $\Hb_{\mathcal{S}}=(\mathcal{S},\mathcal{E}^b_{\mathcal{S}})$ for a set system $\mathcal{S}$ is the hypergraph with vertices identified with the sets in $\mathcal{S}$, and set of edges $\mathcal{E}^{b}_{\mathcal{S}}:=\sett{B\subseteq\mathcal{S}}{B \textrm{ is a minimal blocking set}}$.

			A \emph{coloring} of $\Hb_{\mathcal{S}}$ with at most $q\in \N$ colors is a function $\gamma\ \colon\ \mathcal{S} \rightarrow \range{q}$ such that no edge in $\mathcal{E}^{b}_{\mathcal{S}}$ is monochromatic, i.e., each $B \in \mathcal{E}^b_{\mathcal{S}}$  contains at least two sets of $\mathcal{S}$ colored differently by  $\gamma$.
		\end{definition}
		\begin{proposition}\label{prop:blocking-hg}
		The blocking hypergraph $\Hb_{\mathcal{S}}=(\mathcal{S},\mathcal{E}^b_{\mathcal{S}})$ for a set system $\mathcal{S}$ has the following properties:
			\begin{enumerate}[label=\roman*), ref=\roman*)]
				\item each hyperedge $e$ of $\Hc_\mathcal{S}$ of cardinality at least $3$ induces an edge of $\Hb_{\mathcal{S}}$ of cardinality $2$,\label{prop:blocking-hg:a}
				\item  $\rank{\Hb_{\mathcal{S}}}\leq\rank{\Hc_{\mathcal{S}}}$\label{prop:blocking-hg:b}.
			\end{enumerate}
		\end{proposition}

		Suppose we have a coloring $\gamma\colon\mathcal{S}\to\range{q}$ of the hypergraph $\BB$ with at most $q$ colors.
		For any color $c \in \range{q}$, let $K_c$ be the set of polytopes colored $c$ by $\gamma$, i.e., $K_c = \{ S \in\mathcal{S}\ |\ \gamma(S) = c\}$.
		Since $\gamma$ is a coloring of $\Hb_{\mathcal{S}}$, we immediately have the following.
		\begin{proposition}
			For any color $c \in \range{q}$, $K_c$ does not contain a blocking set as a subset.\label{prop:no_blocking}
		\end{proposition}
		We build a MILP-formulation based on the coloring $\gamma$ of hypergraph $\Hb_{\mathcal{S}}$.
		We introduce a binary variable $z_c\in\left\{0,1\right\}$ for each color class $c$.
		Let $\pi_v$ denote the set of the colors of those polytopes that have $v\in V$ in their vertex set, i.e., 
		$\pi_v = \{ \gamma(P_i)\ |\ v \in S_i \}$.
		Consider the MILP-formulation
		\begin{equation}
			\lambda\in \Delta^V,\quad
				\lambda_v\leq \displaystyle{\sum_{c\in\pi_v}z_c}\ \ \forall v\in V, \quad \displaystyle{\sum_{c=1}^q z_c} = 1,\quad
			z_c \ \in \left\{0,1\right\}\ \ \forall c \in \range{q}.\label{eq:colormip}
		\end{equation}
		A crucial observation about the feasible solutions of (\ref{eq:colormip}) is the following.
		\begin{proposition}\label{prop:cdc-color}
			If $(\lambda,z)$ is a feasible solution of (\ref{eq:colormip}), then $\supp{\lambda}$ does not contain any minimal infeasible set of size at least 3 as a subset.
			\end{proposition}

		\begin{remark}\label{rem:network}
			The LP relaxation of (\ref{eq:colormip}) may have extreme points with fractional  $z_c$ values. However, (\ref{eq:colormip}) admits a network flow representation as defined in \citep{kis2022ideal}, and by~\cite{dobrovoczki_facet_2024}, there exists a polynomial time separation algorithm to identify violated facet defining inequalities.
		\end{remark}
	
	The rank of the blocking hypergraph is at most $d+1$  by Theorem~\ref{thm:conflict-hg} and Proposition~\ref{prop:blocking-hg}. It  can be constructed by enumerating all subsets of simplices of size at most $d+1$, and keeping those that form a minimal blocking set.	
A practical method for coloring the blocking hypergraph is described in Section~\ref{app:coloring} along with an example.
		
		\subsection{Handling pairwise conflicts}\label{sec:pairwise-conflicts}
		Consider the subgraph $G^c_\mathcal{S}$ of rank-two edges of $\Hc_\mathcal{S}$.
		Then the edges of $G^c_\mathcal{S}$ are exactly those pairs of vertices in $V$ where $\lambda$ cannot be positive simultaneously, or in other words, those pairs of vertices that are not contained in any set $S_i$ of $\mathcal{S}$.
		Let $(A_\ell\cup B_\ell, E_\ell)_{\ell=1}^K$ be a biclique cover of $G^c_\mathcal{S}$ and $y_\ell\in\left\{0,1\right\}$ the binary variable corresponding to the choice of biclique $(A_\ell\cup B_\ell, E_\ell)$.
		Then constraints 
		\begin{equation}\label{eq:pairwise-IB}
				\lambda \in \Delta^V,\quad\sum_{v\in A_\ell}\lambda_v \leq  y_\ell,\quad\quad\sum_{v\in B_\ell}\lambda_v \leq  1-y_\ell,\quad\quad y_\ell\in\left\{0,1\right\}\quad\forall \ell\in\range{K}
		\end{equation}
		ensure that for any pair $v_1,\, v_2\in V$ if $v_1,\, v_2\in\supp{\lambda} := \{ v \in V\ |\ \lambda_v > 0\}$, then there exists a set $S_i\in\mathcal{S}$ such that $v_1,\,v_2 \in S_i$.
		Observe that whenever the conflict hypergraph $\Hc_{\mathcal{S}}$ has rank 2, then $G^c_\mathcal{S} = \Hc_\mathcal{S}$ and (\ref{eq:pairwise-IB}) is MILP formulation for (\ref{eq:disj_constr}). The following statement is a direct consequence of the above definitions and the proof is omitted.
		\begin{proposition}\label{prop:cdc-biclique}
			If $(\lambda,y)$ is a feasible solution of (\ref{eq:pairwise-IB}), then $\supp{\lambda}$ does not contain any minimal infeasible set of size 2 as a subset.
		\end{proposition}

	Subsequently, we describe a heuristic procedure to find a biclique cover of a graph $G=(V,E)$ of not too many bicliques.  Our method is based on finding a maximum weight biclique repeatedly. In Algorithm~\ref{alg:biclique_cover}, we start from an empty set $\mathcal{B}=\emptyset$ and the edge-weights are set uniformly to 1, i.e., $w_e=1$ for all $e\in E$.
		We extend $\mathcal{B}$ in each iteration with a new biclique that has the most edges that are not covered by the previous ones, and set $w_e=0$ for edges already covered.
		The new biclique is computed  by solving the mathematical program (\ref{eq:maxBiclique}) with the actual edge weights.
		This is repeated until all edges are covered.
		\begin{algorithm}[!th]
			\small
			\caption{Find biclique cover of graph $G$}\label{alg:biclique_cover}
			\begin{algorithmic}[1]
				\INPUT Graph $G=(V,E)$
				\OUTPUT Biclique cover $\mathcal{B}$
				\State Let $\mathcal{B}=\emptyset$ and $w_e=1$ for all $e\in E$
				\While{$\exists e\in E\colon w_e>0$}
				\State $(A,\,B,\,E') \gets$ maximum weight biclique of $G$ for edge weights $w$
				\State $\mathcal{B}\gets\mathcal{B}\cup\left\{(A,\,B,\,E')\right\}$
				\State $w_e\gets 0$ for all $e\in E'$
				\EndWhile
				\State\Return $\mathcal{B}$
			\end{algorithmic}
		\end{algorithm}
		
	Next we present a MILP formulation for finding a maximum-weight biclique in a graph for arbitrary non-negative edge weights $w\colon E\to \R_{\geq0}$.	
	There are three sets of binary variables, $x^1$, $x^2$, and $y$, where $x^1_u$ and $x^2_u$ indicate whether $u \in V$ belongs to the left or right side of a bipartite subgraph of $G$, and $y_{\{u,v\}}$ indicates if $\{u,v\}$ is an edge of the biclique.
	The formulation is as follows.
		\begin{subequations}
			\begin{align} 
				\min && \sum_{e\in E}w_{e}y_{e}\label{eq:maxBiclique:obj}\\
				\mathrm{s.t.} && x^1_u + x^2_v & \leq  1 &\quad\forall u,v\in V, \{u,v\}\notin E\label{eq:maxBiclique:2}\\
				&& \sum_{u\in V}x^i_u &\geq 1 &\quad\forall i=1,2\label{eq:maxBiclique:3}\\
				&& x^1_u + x^2_v & \leq  1+y_{\{u,v\}} &\quad\forall u,v\in V, \{u,v\} \in E\label{eq:maxBiclique:4}				\\				
				&& x^1_u + x^2_u & \geq  y_{\{u,v\}} &\quad\forall u, v\in V, \{u,v\} \in E\label{eq:maxBiclique:5}\\
				&& x^i_u + x^i_v & \geq  y_{\{u,v\}} &\quad\forall u,v\in V, \{u,v\} \in E,\,\forall i=1,2\label{eq:maxBiclique:6}\\
				&& x^1,\,x^2 &\in \left\{0,1\right\}^V,
				\quad y \in \left\{0,1\right\}^E\label{eq:maxBiclique:8}.
			\end{align}\label{eq:maxBiclique}
		\end{subequations}
		The objective~\eqref{eq:maxBiclique:obj} is to maximize the total weight of those edges covered. 
		For a pair of distinct vertices $u,v$ that do not span an edge, at most one of them can be contained in the biclique~\eqref{eq:maxBiclique:2}, while if $u = v$, then $u$ cannot be on both sides of the biclique.
		Each side of the biclique must have at least one vertex by~\eqref{eq:maxBiclique:3}. 
		If vertices $u$ and $v$ are in different sides of the biclique, then edge $\{u,v\}$ must be an edge of the biclique~\eqref{eq:maxBiclique:4}, while if $y_{\{u,v\}} = 1$, the vertices $u$ and $v$ must be in opposite sides of the biclique by~\eqref{eq:maxBiclique:5} and~\eqref{eq:maxBiclique:6}. The following result follows  from the definitions and the proof is omitted.
		\begin{proposition}
			Every integer solution of~(\ref{eq:maxBiclique}) is a biclique in graph $G$.
		\end{proposition}
		Observe that for weights $w_e=1$ for all $e\in E$, the optimal solution is a maximal biclique in $G$.

		\subsection{Constructing the MILP formulation for (\ref{eq:disj_constr})}
		\label{sec:gib}
		
		The MILP formulation consisting of~(\ref{eq:colormip}) and~(\ref{eq:pairwise-IB}) is called {\em generalized independent branching scheme (GIB)}.
		The next theorem is a direct consequence of Propositions~\ref{prop:cdc-color} and \ref{prop:cdc-biclique}:
		\begin{theorem}
			The GIB formulation (\ref{eq:colormip})-(\ref{eq:pairwise-IB}) is a proper formulation for constraint (\ref{eq:disj_constr}).
		\end{theorem}

		\section{Short-term pumped-storage hydropower plant scheduling and unit commitment}
		\label{sec:hydropower}
		To demonstrate the usefulness of the described techniques, the short-term pumped-storage hydropower plant scheduling and unit commitment problem (STHS) was chosen, with head-dependent hydropower function.
	{\em Pumped-storage hydropower systems\/} are primarily utilized for energy storage and load balancing. They generate power during periods of high demand and electricity prices by discharging water through turbines from a reservoir. Conversely, water is pumped back into the reservoir during low-demand periods when electricity prices are lower.
In short-term hydro scheduling (STHS), the objective is to maximize the profitability of power generation by creating a power generation schedule over a short time frame, typically ranging from a few hours up to a week. This scheduling must comply with various physical, environmental, and regulatory constraints.
The schedule is divided into equal-length time periods. During each time period, three potential actions can be taken: power generation, pumping, or remaining idle. For power generation periods, the water discharge rate must also be determined, while the pumping rate remains fixed.
		
	The {\em hydropower function (HPF)} $\phi$ describes the often non-linear relationship between the water discharge rate and the amount of electric power generated.
	The hydropower function may  also depend on the {\em hydraulic head}, which is defined as the difference in altitude between the upper and lower water levels in the reservoir.
		
		A significant challenge in the modeling of the STHS  is to handle the non-linear and typically non-convex HPF.
		When the hydraulic head is neglected, then the HPF is modeled as a univariate function, see e.g., \citep{chang_experiences_2001, hjelmeland_impact_2018, skjelbred_dynamic_2020}.
		For head dependent HPF 
		quadratic models can be found in~\citep{finardi_solving_2005, catalao_parameterisation_2006}, mixed integer linear programming models are presented by~\citet{borghettiMILPApproachShortTerm2008, britoMixedintegerNonseparablePiecewise2020, alvarezOperationPumpedStorage2020, dossantosabreuContinuousPiecewiseLinear2022}.
		More optimization models for the STHS are collected in a recent survey~\citep{kong_overview_2020}, while different MILP formulations for the HPF are listed and compared by~\cite{guisandez_mixed_2021}.

		Our MILP model strictly follows~\cite{borghettiMILPApproachShortTerm2008}, except for the linearization of $\varphi$, which is done by replacing it with its PWL approximation $\hat{\varphi}$.
		Each $v \in V_\mathcal{T}$ has two coordinates, denoted by $v_q$ and $v_r$, corresponding to the discharge flow rate and the reservoir volume, respectively.
		Then $\hat{\varphi}$ is modeled as follows.
		Let $G$ denote the conflict graph of $\mathcal{T}$, $(A_\ell \cup B_\ell, E_\ell)_{\ell=1}^K$ the biclique cover of $G$ and $\Hb_\mathcal{T}$ the blocking hypergraph of $\mathcal{T}$ as in Definition~\ref{def:blocking}.
		For each time period $t$, let $\lambda^t_v\in\left[0,1\right]$ denote the variable corresponding to the convex coefficient of point $v \in V$, $y^t_\ell\in\left\{0,1\right\}$ the variable for biclique $(A_\ell \cup B_\ell, E_\ell)$ and $z^t_c\in\left\{0,1\right\}$ the variable for color class $c$ in the coloring of $\BB$.
		Let $q_t$ denote the discharge flow rate and $r_t$ the volume of water  in the reservoir in period $t$.
		Denote by $R_{\min}$ and $R_{\max}$ the minimum and maximum reservoir volume.
		The constants $Q^-<0$ and $P^-$ denote the flow pumped by the pump (\SI{}{\cubic\metre\per\second}) and the power consumed by the pump (\SI{}{\mega\watt}) per period, respectively. Let $g_t$ and $u_t$ be the binary variables encoding that the turbine or the pump is used in period $t$, respectively.
		
		The computation of $\varphi(q_t,r_t)=p_{t}$ is modeled as follows:
		\begin{subequations}\label{eq:hpfMILP}
			\begin{equation}
				p_{t}\ =\ \sum_{v\in V}\lambda_v^{t} \varphi(v_q,v_r) + u_{t}P^-\label{constr:hpf:power}
			\end{equation}
			\begin{equation}
			q_{t}\ =\ \sum_{v\in V}\lambda_v^{t} v_q+u_{t}Q^-,\quad R_{\min}u_t\ \leq\ r_t-\sum_{v\in V}\lambda_v^{t} v_r\ \leq\  R_{\max}u_t,\quad\sum_{v\in V}\lambda_v^{t}\ =\ g_{t}\label{constr:hpf:conv}
			\end{equation}
			\begin{equation}
				\sum_{v\in A_i}\lambda_v^{t}\ \leq\ y_i^{t},\quad\sum_{v\in B_i}\lambda_v^{t}\ \leq\ 1-y_i^{t}\quad\forall i\in\range{n}\label{constr:hpf:biclique}
			\end{equation}
			\begin{equation}
				\sum_{c=1}^m z_c^{t}\ =\ g_{t},\quad\sum_{v\in V_\pi} \lambda_v^{t}\ \leq\ \sum_{c\in\pi}z_c^{t}\quad\forall\pi\in\Pi\label{constr:hpf:color}
			\end{equation}
			\begin{equation}
				\lambda_v^{t} \in \left[0,1\right]\quad\forall v\in V,\quad y^{t}_i \in \left\{0,1\right\}\quad\forall i\in\range{n},\quad z^{t}_c \in \left\{0,1\right\}\quad\forall c\in\range{m},\label{constr:hpf:bin}
			\end{equation}
		\end{subequations}
		where $\Pi$ stands for the set of color patterns of the points in $V$, and $V_\pi$ is the set of points with color pattern $\pi\in\Pi$.
		Constraint~\eqref{constr:hpf:power} calculates the generated power and the water discharge as a convex combination in case period $t$ is a generating period and as the power consumed by the pump in case period $t$ is pumping.
		Similarly, constraints~\eqref{constr:hpf:conv} ensure that in case of power generation in period $t$, the reservoir volume and the flow rate is used to calculate the convex coefficients $\lambda$.
		Otherwise, the conservation constraints determine the reservoir volume, while the flow rate is set to a fixed negative pump rate, characteristic to the pump.
		Constraints~\eqref{constr:hpf:biclique} are the pairwise independent branching constraints as in~\eqref{eq:pairwise-IB}, and constraints~\eqref{constr:hpf:color} are the coloring constraints which ensure that there are no conflict sets of cardinality $3$  in $\supp{\lambda}$.
		If the rank of the conflict hypergraph is less then $3$, constraints~\eqref{constr:hpf:color} are not needed, since constraints~\eqref{constr:hpf:biclique} alone resolve all the conflicts.

		\section{Computational experiments}\label{sec:experiment}
		The effectiveness of the developed techniques were assessed through a series of computational experiments conducted on instances of the STHS.
		The non-linear HPF $\varphi$ is a bivariate polynomial, similar to the one approximated in~\cite{borghettiMILPApproachShortTerm2008} and~\cite{thomopulos_generating_2024}:
		\begin{equation}
			\varphi(q,r) = \sum_{h\in H}\left(L_hq\sum_{l\in H} \left(K_lv^l-L_{\mathrm{lb}}-R_0q^2\right)\right),
		\end{equation}
		where $H$ is a set of indices for the coefficients $L_h$ and $K_h$, and $R_0,\, L_{\mathrm{lb}}$ are some constants.
		
		The experiments are divided into two major parts.
		In Section~\ref{sec:experiment:random} we demonstrate the effectiveness of the modeling approach of Section~\ref{sec:hydropower} of $\hat{\varphi}$  by comparing it to different MILP models for $\hat{\varphi}$ based on DLog, Inc, MC, DCC and CC, defined in e.g.~\cite{vielmaMixedIntegerModelsNonseparable2010}. 
		The results presented in Section~\ref{sec:adaptive} are twofold. 
		Firstly, the efficiency of the different PWL function representations are evaluated and compared on the STHS with adaptive triangulations of different accuracy.
		Then the adaptive triangulations are compared to orthogonal grid triangulations of various sizes, as well as the proposed model is compared to the 6-stencil formulation of \cite{huchetteNonconvexPiecewiseLinear2019}.

		The models were implemented in Python, and solved using FICO XPRESS v9.2.5.
		The experiments were performed on 1 thread, with a time limit of 3600 seconds on a server with i9-7960X CPU @ 2.80GHz and Linux operating system.

		\subsection{Random triangulations}\label{sec:experiment:random}

		The models for representing the PWL interpolation of the HPF were tested on three different scenarios of the hydropower scheduling problem, while the PWL interpolation was defined over the random triangulations.
		Each scenario captures a one week long (168, one hour planning periods) horizon of a different month (April, June and December), with different electricity price and inflow predictions and different initial and final reservoir volume requirements.

		We consider two classes of triangulations.
		The first class consists of triangulations with a conflict hypergraph of rank $2$ (called non-blocking triangulations), while the second class contains triangulations with a conflict hypergraph of rank $3$ (referred to as blocking triangulations). Below we present detailed results for non-blocking triangulations only, whereas for blocking triangulations, we refer to Section~\ref{app:blocking}.
		
		For non-blocking triangulations, we use the IB formulation~\eqref{eq:pairwise-IB} based on the biclique covers constructed by the method of Section~\ref{sec:pairwise-conflicts}.
		The sizes of the different triangulations and the corresponding formulations are summarized in Table~\ref{tab:lpsize}.
				
		Table~\ref{tab:unadjusted_correct} shows the runtime of the models for non-blocking triangulation of all sizes for the different scenarios.
		Column \#instances shows the number of triangulations that are non-blocking.
		Columns DLog, Inc, MC, DCC, CC and Our IB show the average runtime of the models on these triangulations, the last two rows showing the average over all problem instances and the number of instances with optimal solution.
		The best runtime is highlighted for each row.
		Our IB formulation dominates on most instance sizes and scenarios, however, in a few cases the other models outperformed it.
		
		\begin{remark}
			Even though the MILP description of the DLog model has seemingly favourable properties, such as the small number of constraints, variables and binaries to the other models (see Table~\ref{tab:lpsize}), it performed poorly in the experiments.
			This is attributed to the large number of non-zeros in the coefficient matrix.
		\end{remark}

		\begin{table}
			\centering
			\small
			\caption{Size of MILP formulation of the models on the random triangulations of different size.\label{tab:lpsize}}
			\footnotesize
				\pgfplotstabletypeset[
				col sep=comma,
				string type,
				columns/size/.style={column name=Size, column type=l},
				columns/category/.style={column name=, column type=r},
				columns/Dlog/.style={column name=DLog, column type=r},
				columns/Inc/.style={column name=Inc, column type=r},
				columns/MC/.style={column name=MC, column type=r},
				columns/DCC/.style={column name=DCC, column type=r},
				columns/CC/.style={column name=CC, column type=r},
				columns/GIB min/.style={column name=min, column type=r},
				columns/GIB max/.style={column name=max, column type=r},
				every head row/.style={
					before row={
						\toprule
						\hline
						&&&&&&& \multicolumn{2}{c@{\hspace{1em}}}{GIB}\\
						\cmidrule{8-9}
					},
					after row=\midrule
				},
				every row no 4/.style={before row=\midrule} ,
				every row no 8/.style={before row=\midrule} ,
				every last row/.style={after row=\bottomrule},
				row sep=\\,
				]{size, category, Dlog, Inc, MC, DCC, CC, GIB min, GIB max\\
					small, rows, 4199,15119, 15623, 6887, 5711, 4871,5879\\
					, columns, 15288,14280, 14448, 18816, 8736, 5544,6048\\
					, binaries, 1512,4872, 5040, 5040, 5040, 1848,2352\\
					, non-zeros, 138596,74252, 74420, 80132, 41492, 34268,39644\\
					medium, rows, 4535,34271, 34775, 13271, 9575, 5543,6887\\
					, columns, 34608,33432, 33600, 44352, 18984, 9744,10416\\
					, binaries, 1680,11256, 11424, 11424, 11424, 2184,2856\\
					, non-zeros, 362708,175556, 175220, 188660, 90212, 76940,87188\\
					large, rows, 4871,65519, 66023, 23687, 15791, 6215,13775\\
					, columns, 66024,64680, 64848, 86016, 35616, 16296,17304\\
					, binaries, 1848,21527, 21840, 21840, 21840, 2520,3528\\
					, non-zeros, 770276,339692, 339860, 365732, 169172, 151700,188156\\
				}
		\end{table}

		\begin{table}[ht]
			\centering
			\caption{Runtime (\SI{}{\second}) of models on the non-blocking triangulations.\label{tab:unadjusted_correct}}
			\footnotesize
				\pgfplotstabletypeset[
				col sep=comma,
				header=true,
				columns/Month/.style={string type, column name=Month, column type=l},
				columns/Size/.style={string type, column name=Size, column type=r},
				columns/inst/.style={string type, column name=\#instances, column type=r},
				columns/Dlog/.style={column name=DLog, fixed, precision=2, column type=r, fixed zerofill=true},
				columns/Inc/.style={column name=Inc, fixed, precision=2, column type=r, fixed zerofill=true},
				columns/MC/.style={column name=MC, fixed, precision=2, column type=r, fixed zerofill=true},
				columns/DCC/.style={column name=DCC, fixed, precision=2, column type=r, fixed zerofill=true},
				columns/CC/.style={column name=CC, fixed, precision=2, column type=r, fixed zerofill=true},
				columns/GIB/.style={column name=Our IB, fixed, precision=2, column type=r, fixed zerofill=true},
				every head row/.style={before row={\toprule\hline}, after row=\midrule},
				every row no 3/.style={before row=\midrule},
				every row no 6/.style={before row=\midrule},
				every row no 9/.style={before row=\midrule},
				every row no 10/.style={before row=\midrule},
				every last row/.style={after row=\bottomrule},
				every row 0 column 8/.style={postproc cell content/.style={@cell content/.add={$\bf}{$}}},
				every row 1 column 8/.style={postproc cell content/.style={@cell content/.add={$\bf}{$}}},
				every row 2 column 8/.style={postproc cell content/.style={@cell content/.add={$\bf}{$}}},
				every row 3 column 7/.style={postproc cell content/.style={@cell content/.add={$\bf}{$}}},
				every row 4 column 8/.style={postproc cell content/.style={@cell content/.add={$\bf}{$}}},
				every row 5 column 8/.style={postproc cell content/.style={@cell content/.add={$\bf}{$}}},
				every row 6 column 6/.style={postproc cell content/.style={@cell content/.add={$\bf}{$}}},
				every row 7 column 8/.style={postproc cell content/.style={@cell content/.add={$\bf}{$}}},
				every row 8 column 8/.style={postproc cell content/.style={@cell content/.add={$\bf}{$}}},
				every row 9 column 8/.style={postproc cell content/.style={@cell content/.add={$\bf}{$}}},
				opt/.list={3,4,5,6,7},
				boldopt/.list={8},
				row sep=\\,
				]{
					Month,Size,inst,Dlog,Inc,MC,DCC,CC,GIB\\
					April,small,8,1889.20375,831.3025,185.78625,63.06,30.6775,10.40875\\
					,medium,9,3600,3401.587778,3267.922222,1879.234444,1728.998889,542.8322222\\
					,large,4,3600,3600,3600,3543.0525,3600,1152.0225\\
					June,small,8,940.26875,459.9825,65.46,474.81625,32.02375,89.51125\\
					,medium,9,1363.997778,1171.035556,1028.356667,865.6344444,814.11,434.4588889\\
					,large,4,2714.12,1501.6025,2807.5275,1450.3575,333.71,29.955\\
					December,small,8,1836.62625,592.615,733.225,280.0725,449.61125,457.01\\
					,medium,9,3600,2913.986667,2986.436667,2950.347778,2653.421111,2159.376667\\
					,large,4,3600,2726.505,3600,2878.545,1955.18,226.285\\
					Average time,,,2445.530476,1805.765,1800.860159,1417.275556,1181.324444,608.2301587\\
					\# Optimal,,63,21,36,35,42,47,55\\
				}
		\end{table}

		\subsection{Adaptive triangulations}\label{sec:adaptive}
		
		Adaptive triangulations of different sizes were generated for the HPF, as described in Section~\ref{sec:pwl_fit}.
		Here $\Omega=\left[Q_{\min},\,Q_{\max}\right]\times \left[R_{\min},\,R_{\max}\right]$, $V_0 = \left\{Q_{\min},\,Q_{\max}\right\}\times\left\{R_{\min},\,R_{\max}\right\}$ and $\mathcal{T}_0$ is one of the two possible triangulations of $\Omega$ with $V(\mathcal{T}_0)=V_0$.
		Six different error tolerance values were used: $\varepsilon\in\left\{0.5,0.4,0.3,0.2,0.1,0.05\right\}$.
		
		The sizes of the triangulations corresponding to the error tolerances (given unique IDs {\ttfamily a05}, {\ttfamily a04}, {\ttfamily a05}, {\ttfamily a02}, {\ttfamily a01} and {\ttfamily a005}, respectively) are summarized in Table~\ref{tab:adaptiveSize}.
		Figure~\ref{fig:maxError} shows the convergence of the maximum empirical error on a logarithmic scale.
		The non-monotonicity and the observed peaks are due to that the subsequent triangulations are not refinements of each other.
		The resulting triangulation is shown in Figure~\ref{fig:a005}.
		Notably, none of the generated adaptive triangulation instances contained blocking sets of triangles, hence, the methods described in Section~\ref{sec:high-rank-conflict} were not applied to them.

		For comparison, orthogonal grid triangulations were generated.
		Both axes are subdivided into $4,\,8,\,16$ and	$32$ parts by $5,\,9,\,17$ and $33$ equidistant points, named {\ttfamily g4}, {\ttfamily g8}, {\ttfamily g16} and {\ttfamily g32}, respectively.
		The subdivision of grid cells into triangles are chosen randomly.
		Suppose that $d_1,d_2$ is the number of breakpoints on the two axes that define a grid, $\mathcal{T}$ is a grid triangulation.
		It is shown by~\cite{huchetteNonconvexPiecewiseLinear2019} that the conflict graph of $\mathcal{T}$ can be covered by $\log_2(d_1)+\log_2(d_2)+6$ bicliques and the authors gave a construction for such a cover, known as the 6-stencil formulation.
		Triangulation sizes are collected in Table~\ref{tab:adaptiveSize} for both adaptive and grid triangulations, as well as the size of the biclique cover for the adaptive triangulations and the maximal empirical error.
		For the adaptive triangulations, this upper bound is the error tolerance used in their construction, while for the grid triangulations the maximal error is estimated similarly to the method used by Algorithm~\ref{alg:fitting}.

		The runtime of the models on the three scenarios and six adaptive triangulations are shown in Table~\ref{tab:adaptiveCompare}.
		Similarly to the experiment on random triangulations, our IB formulation dominates the others, in this case on all problem instances.

		The quality of the adaptive triangulations were tested against the grid triangulations, as well as our IB formulation against the 6-stencil formulation.
		Table~\ref{tab:adaptiveSize} shows that the 6-stencil formulation and our IB formulation result in similar size MILP models for triangulations of similar sizes.
		In table~\ref{tab:adaptiveResults} we compare the performance of the methods.
		Column PWL obj. contains the optimal objective value of the MILP model,
		while column NL obj. contains the value of the objective function with original HPF evaluated on the optimal values of $q$ and $r$ variables.
		Column Rel.err. is the relative error between the PWL objective value and the non-linear objective value.
		It shows that as the triangulations get finer, the PWL objective converges to the non-linear objective.
		Column Avg.abs. HPF err. shows the absolute deviation between the PWL interpolation of the HPF and the non-linear HPF averaged in the generating periods.
		The runtime and non-linear objective for largest grid and adaptive triangulations ({\ttfamily g32} and {\ttfamily a005}) are highlighted, to show that a very similar non-linear objective value can be achieved by the two models, while the runtime of our IB is much lower.
		\begin{table}[h!]
			\centering
			\caption{Size and measured maximal absolute error of adaptive and grid triangulations and size of the resulting MILP formulations.\label{tab:adaptiveSize}}
			\footnotesize
				\pgfplotstabletypeset[
				col sep=comma,
				columns/Type/.style={string type}, 
				columns/ID/.style={string type}, 
				columns={Type, ID, GridSize, Ntrg, Npts, Nbiclique, AbsErr,Rows,Cols,Bins,NZs}, 
				every head row/.style={before row=\toprule, after row=\midrule}, 
				every last row/.style={after row=\bottomrule}, 
				columns/Type/.style={string type, column type=l}, 
				columns/ID/.style={string type, column type=r}, 
				columns/GridSize/.style={string type, column type=c, column name=Grid size}, 
				columns/Ntrg/.style={column type=r, column name=\#triangles}, 
				columns/Npts/.style={column type=r, column name=\#points}, 
				columns/Nbiclique/.style={string type, column type=c,column name=\#bicliques}, 
				columns/AbsErr/.style={column type=l, precision=4, fixed, fixed zerofill=true, column name=Abs. error, postproc cell content/.style={@cell content/.add={$\leq$ }{}}}, 
				every head row/.style={
					before row={
						\toprule
						\hline
					},
					after row=\midrule
				},
				columns/Rows/.style={fixed, column name=Rows, column type=r},
				columns/Cols/.style={fixed, column name=Columns, column type=r},
				columns/Bins/.style={fixed, column name=Binaries, column type=r},
				columns/NZs/.style={fixed, column name=Non-zeros, column type=r},
				every row no 4/.style={before row=\midrule}, 
				every head row/.style={
					before row={
						\toprule
						\hline
						&&&\multicolumn{4}{c@{\hspace{1em}}}{Triangulation} &\multicolumn{4}{c@{\hspace{1em}}}{MILP size}\\
						\cmidrule(lr){4-7}\cmidrule(lr){8-11}
					},
					after row=\midrule
				},
				row sep=\\,
				]{
					Type, ID, GridSize, Ntrg, Npts, Nbiclique, AbsErr,Rows,Cols,Bins,NZs\\
					grid, {\ttfamily g4}, $4\times 4$, 32, 25, -, 0.9263,5879,7224,2352,44516\\
					, {\ttfamily g8}, $8\times 8$, 128, 81, -, 0.3736,6551,16968,2688,152204\\
					, {\ttfamily g16}, $16\times 16$, 512, 289, -, 0.1628,7223,52248,3024,613868\\
					, {\ttfamily g32}, $32\times 32$, 2048, 1089, -, 0.0366,7895,186984,3360,2656244\\
					adaptive, {\ttfamily a05}, -, 26, 21, 7, 0.5,4871,6048,1848,38804\\
					, {\ttfamily a04}, -, 33, 25, 8, 0.4,5207,6888,2016,49724\\
					, {\ttfamily a03}, -, 42, 30, 9, 0.3,5543,7896,2184,60476\\
					, {\ttfamily a02}, -, 58, 40, 9, 0.2,5543,9576,2184,74420\\
					, {\ttfamily a01}, -, 137, 83, 11, 0.1,6215,17136,2520,160772\\
					, {\ttfamily a005}, -, 306, 176, 13, 0.05,6887,33096,2856,382700\\
				}
		\end{table}

		\begin{table}[h!]
			\centering
			\small
			\caption{Runtime (s) of models on the adaptive triangulations.\label{tab:adaptiveCompare}}
			\footnotesize
				\pgfplotstabletypeset[
				col sep=comma,
				header=true,
				columns/Month/.style={string type, column name=Month, column type=l},
				columns/ID/.style={string type, column name=ID, column type=r},
				columns/Dlog/.style={column name=DLog, fixed, precision=2, column type=r, fixed zerofill=true},
				columns/Inc/.style={column name=Inc, fixed, precision=2, column type=r, fixed zerofill=true},
				columns/MC/.style={column name=MC, fixed, precision=2, column type=r, fixed zerofill=true},
				columns/DCC/.style={column name=DCC, fixed, precision=2, column type=r, fixed zerofill=true},
				columns/CC/.style={column name=CC, fixed, precision=2, column type=r, fixed zerofill=true},
				columns/GIB/.style={column name=Our IB, fixed, precision=2, column type=r, fixed zerofill=true},
				every head row/.style={before row={\toprule\hline}, after row=\midrule},
				every row no 6/.style={before row=\midrule},
				every row no 12/.style={before row=\midrule},
				every row no 18/.style={before row=\midrule},
				every last row/.style={after row=\bottomrule},
				row sep=\\,
				]{
					Month,ID,Dlog,Inc,MC,DCC,CC,GIB\\
					April,{\ttfamily a05},3600,550.22,132.09,59.62,25.88,7.66\\
					,{\ttfamily a04},1106.94,3600,412.34,22.7,57.12,9.44\\
					,{\ttfamily a03},521.93,190.37,1256.79,60.49,176.8,16.33\\
					,{\ttfamily a02},3600,3600,253.99,96.82,57.19,8.49\\
					,{\ttfamily a01},1951.69,3600,1248.7,386.4,498.32,12.62\\
					,{\ttfamily a005},3600,3600,3600,1892.07,2462.73,35.72\\
					June,{\ttfamily a05},23.26,5.67,38.05,7.67,6.97,2.97\\
					,{\ttfamily a04},92.09,2.93,41.84,12.01,8.13,1.19\\
					,{\ttfamily a03},207.57,14.89,66.5,17.8,12.23,2.12\\
					,{\ttfamily a02},3600,9.34,109.03,107.11,22.21,2.5\\
					,{\ttfamily a01},21.75,64.19,263.66,98.71,41.51,12.5\\
					,{\ttfamily a005},251.78,277.03,285.15,264.14,56.9,19.7\\
					December,{\ttfamily a05},62.04,81.2,69.45,6.53,5.65,3.68\\
					,{\ttfamily a04},19.8,22.94,32.9,5.4,3.24,2.39\\
					,{\ttfamily a03},30.08,24.82,16.56,5.84,3.27,1.35\\
					,{\ttfamily a02},3600,3600,136.08,37.14,25.56,7.02\\
					,{\ttfamily a01},315.05,364.93,173.39,103.61,36.95,1.63\\
					,{\ttfamily a005},192.99,3600,763.72,252.67,20.69,2.25\\
					Average,,1266.541111,1289.81,494.4883333,190.9294444,195.6305556,8.308888889\\
				}
		\end{table}

		\begin{table}[h!]
			\centering
			\caption{Results of short-term hydropower scheduling problem of a pumped storage hydropower plan with PWL interpolation of the head dependent HPF. The HPF is modeled with the our IB and 6-stencil models on adaptive and grid triangulations (respectively).\label{tab:adaptiveResults}}
			\footnotesize
				\pgfplotstabletypeset[
				col sep=comma,
				columns/Month/.style={
					column type=l,
					string type
				}, 
				columns/Model/.style={
					column type=r,
					string type
				}, 
				columns/Type/.style={
					column type=r,
					string type
				}, 
				columns/ID/.style={
					string type,
					column type={>{\ttfamily }r},
					column name=\normalfont{ID}
				},  
				columns/Runtime/.style={
					column type=r,
					column name=Runtime (s)}, 
				columns/Obj/.style={
					column type=r,
					column name=PWL obj.,
					fixed, precision=2},     
				columns/NLObj/.style={
					column type=r, 
					column name=NL obj.,
					fixed, precision=2},   
				columns/relerr/.style={
					column type=r,
					column name=Rel.err.,
					dec sep align,
					preproc/expr={100*##1},
					postproc cell content/.append code={
						\ifnum1=\pgfplotstablepartno
						\pgfkeysalso{@cell content/.add={}{\%}}%
						\fi
					},
					precision = 2,
					fixed,
					fixed zerofill=true
				}, 
				columns/avghpferr/.style={
					column type=r,
					fixed,
					fixed zerofill=true,
					column name={\makecell[b]{Avg. abs.\\ HPF err.}},
					precision=4}, 
				every head row/.style={before row=\toprule, after row=\midrule}, 
				every last row/.style={after row=\bottomrule}, 
				columns={Month, Model, Type, ID, Runtime, Obj, NLObj, relerr, avghpferr}, 
				every head row/.style={
					before row={
						\toprule
						\hline
					},
					after row=\midrule,
					column type={>{\rmfamily}r}
				},
				every head row column 3/.style={column type={>{\rmfamily}r}},
				every row no 3/.style={after row={\cmidrule{2-10}}},
				every row no 13/.style={after row={\cmidrule{2-10}}},
				every row no 23/.style={after row={\cmidrule{2-10}}},
				every row no 10/.style={before row=\midrule},
				every row no 20/.style={before row=\midrule},
				every row 29 column 6/.style={
					postproc cell content/.style={
						@cell content/.add={$\bf}{$}
					}
				},
				every row 23 column 6/.style={
					postproc cell content/.style={
						@cell content/.add={$\bf}{$}
					}
				},
				every row 19 column 6/.style={
					postproc cell content/.style={
						@cell content/.add={$\bf}{$}
					}
				},
				every row 13 column 6/.style={
					postproc cell content/.style={
						@cell content/.add={$\bf}{$}
					}
				},
				every row 9 column 6/.style={
					postproc cell content/.style={
						@cell content/.add={$\bf}{$}
					}
				},
				every row 3 column 6/.style={
					postproc cell content/.style={
						@cell content/.add={$\bf}{$}
					}
				},
				every row 29 column 4/.style={
					postproc cell content/.style={
						@cell content/.add={$\bf}{$}
					}
				},
				every row 23 column 4/.style={
					postproc cell content/.style={
						@cell content/.add={$\bf}{$}
					}
				},
				every row 19 column 4/.style={
					postproc cell content/.style={
						@cell content/.add={$\bf}{$}
					}
				},
				every row 13 column 4/.style={
					postproc cell content/.style={
						@cell content/.add={$\bf}{$}
					}
				},
				every row 9 column 4/.style={
					postproc cell content/.style={
						@cell content/.add={$\bf}{$}
					}
				},
				every row 3 column 4/.style={
					postproc cell content/.style={
						@cell content/.add={$\bf}{$}
					}
				},
				row sep=\\,
				]{
					Month, Model, Type, ID, Runtime, Obj, NLObj, relerr, avghpferr\\
					April, 6-stencil, grid, g4, 7.12, 37841.51778, 38038.35694, 0.005175, 0.0131\\
					,,, g8, 9.52, 38846.86542, 38903.96587, 0.001468, 0.0038\\
					,,, g16, 44.46, 39402.04492, 39421.24091, 0.000487, 0.0012\\
					,,, g32, 129.17, 39561.21818, 39564.78645, 0.000090, 0.0002\\
					,Our IB,adaptive, a05, 7.66, 38618.57335, 39438.30104, 0.020785, 0.0551\\
					,,, a04, 9.44, 38777.09938, 39378.65219, 0.015276, 0.0411\\
					,,, a03, 10.11, 38843.45231, 39403.17621, 0.014205, 0.0369\\
					,,, a02, 8.49, 39232.75134, 39473.24995, 0.006093, 0.0148\\
					,,, a01, 12.62, 39464.51399, 39530.59739, 0.001672, 0.0046\\
					,,, a005, 35.72, 39512.77692, 39556.72837, 0.001111, 0.0029\\
					June, 6-stencil,grid, g4, 0.55, 70916.95163, 71514.08472, 0.008350, 0.04\\
					,,, g8, 1.9, 72184.31151, 72438.43967, 0.003508, 0.0155\\
					,,, g16, 8.77, 73164.8184, 73207.38746, 0.000581, 0.0028\\
					,,, g32, 38.28, 73338.75071, 73349.29195, 0.000144, 0.0007\\
					,Our IB,adaptive, a05, 2.97, 72203.4773, 72935.07805, 0.010031, 0.0484\\
					,,, a04, 1.19, 72203.28582, 72924.33839, 0.009888, 0.0476\\
					,,, a03, 2.12, 72715.05697, 73205.65273, 0.006702, 0.0315\\
					,,, a02, 2.5, 72715.57007, 73204.6097, 0.006680, 0.0316\\
					,,, a01, 12.5, 73151.94555, 73300.2732, 0.002024, 0.0099\\
					,,, a005, 19.7, 73276.18162, 73343.79155, 0.000922, 0.0045\\
					December, 6-stencil,grid, g4, 3.22, 173692.0656, 174305.492, 0.003519, 0.0311\\
					,,, g8, 0.57, 176428.6545, 176737.1803, 0.001746, 0.0186\\
					,,, g16, 2.14, 177713.3605, 177743.5356, 0.000170, 0.0018\\
					,,, g32, 14.72, 177893.1674, 177906.7622, 0.000076, 0.0006\\
					,Our IB,adaptive, a05, 3.68, 175829.2586, 177211.5298, 0.007800, 0.0749\\
					,,, a04, 2.39, 176679.438, 177709.2751, 0.005795, 0.0630\\
					,,, a03, 1.35, 176915.1831, 177699.2411, 0.004412, 0.0472\\
					,,, a02, 7.02, 177230.6748, 177712.4175, 0.002711, 0.0263\\
					,,, a01, 1.63, 177648.1492, 177849.7648, 0.001134, 0.0103\\
					,,, a005, 2.25, 177838.5122, 177928.7203, 0.000507, 0.0053\\
				}
		\end{table}
		
		\begin{figure}[h]
			\footnotesize
			\centering
			\subcaptionbox{\footnotesize Convergence of measured error.\label{fig:maxError}}
			{

			}
			\hfill
			\subcaptionbox{\footnotesize Adaptive triangulation {\ttfamily a005}.\label{fig:a005}}
			{
			\centering
			%
		}
		\caption{Convergence of the maximal absolute error of the PWL interpolation of the HPF in the process of the construction of triangulation {\ttfamily a005}.}
		\end{figure}

		\begin{figure}[h]
			\centering
			%
			\caption{Optimal schedule and unit commitment over triangulation {\ttfamily{a005}} for June.\label{fig:June}}
		\end{figure}

		
		%
		\section{Conclusions}
		
		In this paper, we presented computational methods for constructing a PWL approximation of a
		multivariate non-linear function as well as a MILP formulation for continuous PWL functions over
		simplicial partitions of bounded domains in $\mathbb{R}^d$ for arbitrary dimension $d$. Our method for PWL
		function fitting was based on Lipschitz correction and sampling the domain for estimating the approximation error. To reduce the complexity of the MILP formulation,we devised an original rank
		reduction method for the underlying conflict hypergraph and algorithms for finding biclique covers
		of ordinary graphs. Finally,we assessed the performance of the proposed techniques through a series
		of computational experiments on the STHS problem, where the non-linear bivariate hydropower
		function was represented by our model. Through the experiments, we found that our formulation
		outperforms the known applicable models. The computational experiments on the biclique covers
		suggest that co-planar graphs might be covered with a number of bicliques that is logarithmic in
		the number of edges; however, there is no theoretical guarantee for that yet. Further research on
		Lipschitz-correction algorithms for particular families of functions is subject to future work.
		
		\newpage
		\section*{Appendix}
		\appendix

		\numberwithin{proposition}{section}
		\numberwithin{example}{section}
		\numberwithin{algorithm}{section}
		\numberwithin{table}{section}
		\numberwithin{figure}{section}
		\section{Supplementary materials for Section~\ref{sec:pwl_fit}}		
		\label{app:fitting_proofs}
		\subsection{Proofs of Theorems~\ref{thm:lipschitz-bound}, \ref{thm:max-lipschitz}, and \ref{thm:algorithm}}
		\label{app:lipschitz_corr}
			For proving Theorem~\ref{thm:lipschitz-bound}, we need the following result.
			\begin{proposition}\label{prop:simplex-faces}
				For any $p\in\R^d$, there exist faces $F_1,\,F_2$ of $T$
				and points $x\in F_1,\, y\in F_2$ such that
				$F_1\cap F_2 =\emptyset$, and  $p$ is parallel to $y-x$.
			\end{proposition}
			\begin{proof}
				Since $T$ is full dimensional, $V = V(T)$ has cardinality $d+1$.
				Consider the following linear program:
				\begin{subequations}
					\begin{equation}
						\max\ \gamma \label{eq:simplex-faces:obj}\\
					\end{equation}
					\begin{equation}
						\mathrm{s.t.}\quad y = x + \gamma v\label{eq:simplex-faces:1}
					\end{equation}
					\begin{equation}
						x = \sum_{u\in V}\lambda_uu,\,\sum_{u\in V}\lambda_u = 1,\,\lambda_u \geq 0\ \forall u\in V\label{eq:simplex-faces:2}
					\end{equation}
					\begin{equation}
						y = \sum_{w\in V}\mu_ww,\,\sum_{w\in V}\mu_w = 1,\,\mu_w \geq 0\ \forall w\in V\label{eq:simplex-faces:3}
					\end{equation}
				\end{subequations}
				Equations~(\ref{eq:simplex-faces:2}) and~(\ref{eq:simplex-faces:3}) express that $x$ and $y$ lie in $T$, while equation~(\ref{eq:simplex-faces:1}) together with the objective function ensures that $y-x$ is parallel to $p$. Clearly, any vertex $u$ of $T$ induces  a feasible solution $(x=u, y=u, \gamma = 0, \lambda = e_u, \mu = e_u)$ of this LP.
				There are $4d+3$ variables, $2d+2$ of them non-negative, and $3d+2$ equality constraints.
				Hence, in a basic solution $d+1$ of the non-negative variables equal to $0$, and at most $d+1$ of them can be positive.
				Let $\lambda^*,\,\mu^*$ denote the optimal solution for $\lambda$ and $\mu$, respectively, and let $V_1=\supp{\lambda^*}$ and $V_2=\supp{\mu^*}$.
				If $V_1\cap V_2\neq\emptyset$, then there is a vertex of $T$ that is not used in the convex combination for  $x$ nor for $y$.
				Hence, $x$ and $y$ lie on a face $F$ of $T$, and $p$ is parallel to $F$, in which case we repeat the argument using $F$ in place of $T$.
				Otherwise, we have that $V_1$ and $V_2$ span two disjoint faces of $T$, and they contain all $d+1$ vertices of it. Moreover, $x$ is in the face $F_1$ spanned by $V_1$ and $y$ is in the face $F_2$ spanned by $V_2$, as claimed.
			\end{proof}
	
			\begin{proof}[Proof of Theorem~\ref{thm:lipschitz-bound}]
				Since $\hat{f}$ is an affine function on each $T \in \mathcal{T}$, there exists a number $b_T \in \mathbb{R}$ such that $\hat{f}_T(x) = \nabla\hat{f}_T \cdot x + b_T$ for all $x \in T$.
				Hence, the directional derivative of $\hat{f}_T$ is  $D_v \hat{f}(x) = \nabla\hat{f}_T \cdot v$ for all $x \in T$ and unit vector $v \in \mathbb{R}^d$.  
				Now we apply Proposition~\ref{prop:simplex-faces} to the vector $p = \nabla \hat{f}_T / \lVert \nabla \hat{f}_T \rVert$ to obtain two vectors $x,y \in T$ such that $y-x$ is parallel to $p$ and $x$, $y$ lie in disjoint faces $F_x$ and $F_y$ of $T$. 
				By the above definitions, we have 
				$\lVert \nabla \hat{f}_T \rVert = \nabla \hat{f}_T \cdot p = |\nabla \hat{f}_T \cdot (y-x)| / \lVert y-x \rVert = |\hat{f}_T(y) - \hat{f}_T(x)| / \lVert y-x \rVert$. 
				
				Let $V_x$ and $V_y$ be the set of vertices of the faces $F_x$ and $F_y$, respectively. Since $x$ and $y$ are convex combinations of the vertices of $F_x$ and $F_y$, respectively, and $f(u) = \hat{f}(u)$ for all $u \in V_x \cup V_y$, we have
				\[
				|\hat{f}_T(y)-\hat{f}_T(x)|  \leq \max_{u\in V_x,\, w\in V_y}\lvert f(u)-f(w)\rvert
				\leq  \max_{u\in V_x,\, w\in V_y} L \cdot \lVert u-w\rVert\leq L\cdot\ell_T^{\max}.
				\]
				On the other hand, $\lVert y-x \rVert \geq \delta^{\min}_T$ by the definition of $\delta^{\min}_T$.
				Hence, we have $\lVert\nabla\hat{f}_T\rVert\leq L\cdot\nicefrac{\ell^{\max}_T}{\delta^{\min}_T}$.
			\end{proof}

		\begin{proof}[Proof of Theorem~\ref{thm:max-lipschitz}.]
			As in the proof of Theorem~\ref{thm:lipschitz-bound},
			we can prove that there exist $x,y \in T$ that lie in disjoint faces of $T$ and 	$\lVert \nabla \hat{f}_T \rVert = |\hat{f}_T(y) - \hat{f}_T(x)| / \lVert y-x \rVert$.	
			Since $T$ is a triangle in $\mathbb{R}^2$, $x$ or $y$ is a vertex of $T$.
			So, suppose $y$ is a vertex of $T$, and then $x$ is a point on the opposite edge of $T$.
			Let $u$ and $w$ be the other two vertices of $T$.
			Denote $\ell_y$ the length of the edge opposite to $y$, $h_y$ and $\alpha_y$ the altitude and angle corresponding to $y$. We have similar notation for $u$ and $w$ as well.
			Recall that the area of $T$ is $\ell_u h_u/2=\ell_y h_y / 2 =\ell_w h_w/2$.
			Using $h_y = \ell_w \sin \alpha_u = \ell_u \sin \alpha_w$, we derive $\ell_u / h_y =  1/\sin \alpha_w$
			and $\ell_w / h_y = 1 / \sin \alpha_u$.  Now we compute
			\[			
			\begin{split}
				\lVert \nabla \hat{f}_T \rVert & = \frac{|\hat{f}(y)-\hat{f}(x)|}{\lVert y-x\rVert} \overset{(1)}\leq \frac{\max\left\{|f(y)-f(u)|,|f(y)-f(w)|\right\}}{\lVert y-x\rVert}	\overset{(2)}\leq \frac{\max\left\{L\lVert y-u\rVert,L\lVert y-w\rVert\right\}}{\lVert y-x\rVert}\\
				& \overset{(3)}\leq  L\cdot\frac{\max\left\{\ell_u, \ell_w\right\}}{h_y}
				\overset{(4)}= L\cdot\max\left\{\frac{1}{\sin\left(\alpha_u\right)},\frac{1}{\sin\left(\alpha_w\right)}\right\}\leq \frac{L}{\sin\left(\alpha_T^{\min}\right)},
			\end{split}
			\]
			where (1) follows from the fact that $\hat{f}(x)$ is the convex combination of $f(u)$ and $f(w)$ and $f$ equals $\hat{f}$ on the vertices of $T$.
			The Lipschitz-continuity of $f$ implies (2).  
			Since the altitude $h_y$ is the minimum  distance of $y$ to the line passing through $u$ and $w$, we have $\lVert y-x\rVert\geq h_y$. Moreover, $\ell_w = \lVert y-u\rVert$ and $\ell_u = \lVert y-w\rVert$, and  inequality (3) holds.
			Equation (4) follows from our expressions for $\ell_u / h_y$ and $\ell_w / h_y$.
			Hence, $\lVert\nabla\hat{f}_T\rVert \leq\nicefrac{L}{\sin\left(\alpha^{\min}_T\right)}$.
		\end{proof}

			\subsection{Pseudocode for sub-routines of Algorithm~\ref{alg:fitting}}\label{app:pseudocodes_subroutines}
			In this section we provide high-level pseudocode for the Delaunay refinement and sampling sub-routines that appear in Algorithm~\ref{alg:fitting}. For the detailed description of Ruppert's algorithm for Delaunay refinement refer to~\cite{rupert1995delaunay} and for the maximal Poisson-disk sampling refer to~\cite{ebeida_simple_2012}.
			
			\begin{minipage}[t]{0.44\textwidth}
				\renewcommand{\thealgorithm}{Refinement}
				\begin{algorithm}[H]
					\small
					\caption{Ruppert's Delaunay refinement algorithm~\citep{rupert1995delaunay}}\label{app:ruppert_alg}
					\begin{algorithmic}[1]
						\INPUT Delaunay triangulation $\mathcal{T}$ of $\Omega\subset\R^2$, angle bound $0<\alpha<20$.
						\OUTPUT Delaunay triangulation $\mathcal{T}$ with smallest angle at least $\alpha$.

						\Repeat
						\While{any edge $e$ encroached upon}
						\State Split edge $e$, update $\mathcal{T}$
						\EndWhile
						\State Let $T \in \mathcal{T}$ be a skinny triangle (min angle $< \alpha$)
						\State $p \gets$ circumcenter of $T$
						\If{$p$ encroaches upon any edges $e_1,\ldots,e_k$}
						\State Split edges $e_1,\ldots,e_k$,  update $\mathcal{T}$
						\Else
						\State Split triangle $T$ with $p$, update $\mathcal{T}$
						\EndIf
						\Until{no edges encroached upon and no angles $\leq \alpha$}
						
						\State \Return $\mathcal{T}$
					\end{algorithmic}
				\end{algorithm}
				
			\end{minipage}
			\hfill
			\begin{minipage}[t]{0.5\textwidth}
				\renewcommand{\thealgorithm}{MPS}
				\begin{algorithm}[H]
					\small
					\caption{Maximal Poisson-disk sampling~\citep{ebeida_simple_2012}.\label{alg:mps}}
					\begin{algorithmic}[1]
						\INPUT Simplex $T\subset\R^d$, covering radius $r_T>0$ and minimal distance of sample points $r_{\min}$.
						\OUTPUT Set of sample points $X_T$ such that $X_T$ covers $T$ with radius $r_T$ with $\min_{x,y\in X_T}\lVert x-y\rVert > r_{\min}$
						\State Construct background grid $\mathcal{G}$ with cell diameter $h \le r_T/\sqrt{d}$
						\State $X_T \gets \emptyset$
						\State $\mathcal{A} \gets$ set of all grid cells intersecting $T$ \Comment{Active cells}			
						\While{$\mathcal{A} \neq \emptyset$}\Comment{Start covering active cells}
						\ForAll{active cells $C \in \mathcal{A}$}
						\State Generate a sample point $p \in C \cap T$
						\If{$\min_{x\in X_T}\lVert x-p\rVert > r_{\min}$}
						\State $X_T\gets X_T\cup\{p\}$
						\State $\mathcal{A}\gets \mathcal{A}\setminus\{C\}$
						\EndIf
						\EndFor
						\State $\mathcal{A}'\gets \emptyset$\Comment{Start refinement}
						\ForAll{$C \in \mathcal{A}$}
						\State Refine $C$ into $2^d$ subcells
						\ForAll{subcells $C'$}
						\If{$C'$ is not fully covered by $X_T$}
						\State Add $C'$ to $\mathcal{A}'$
						\EndIf
						\EndFor
						\EndFor
						\State $\mathcal{A} \gets \mathcal{A}'$
						\EndWhile
						\State\Return $X_T$
					\end{algorithmic} 
				\end{algorithm}
			\end{minipage}
			
		\subsection{Proof of Proposition~\ref{prop:sample-error} and Theorem~\ref{thm:algorithm}}\label{app:sampling}
		
		\begin{proof}[Proof of Proposition~\ref{prop:sample-error}.] 
			For any $p \in T$, let $\sigma_T(p)=\min_{x\in X_T}\lVert x-p \rVert$ and $\hat{e}_T = \max_{x \in X_T} |f(x) - \hat{f}_T(x)|$.
			Since $X_T$ covers $T$ with radius $r_T$, we have $\sigma_T(p) \leq r_T$ for any $p \in T$.
			Let $p^*$ be the maximizer of $\max_{p\in T}|f(p)-\hat{f}(p)|$, and $x\in X_T$  the sample point such that $\sigma_T(p^*)= \lVert x-p^*\rVert$. By the definition of $r_T$, we have $\sigma_T(p^*)\leq r_T$. Then we compute
			\begin{equation*}
				\begin{aligned}
					\varepsilon_T = |f(p^*)-\hat{f}(p^*)|& = |f(p^*)-f(x)+f(x)-\hat{f}(x)+\hat{f}(x)-\hat{f}(p^*)|\\
					&\leq |f(p^*)-f(x)|+|f(x)-\hat{f}(x)|+|\hat{f}(x)-\hat{f}(p^*)|\\
					&\overset{(1)}\leq \sigma_T(p^*)L+\hat{\varepsilon}_T+\sigma_T(p^*)\hat{L}_T\ \leq\  \hat{\varepsilon}_T+(L+\hat{L}_T)r_T,
				\end{aligned}
			\end{equation*}
			where (1) follows from the Lipschitz-continuity of $f$ and $\hat{f}$, and from the definition of $\hat{e}_T$.
		\end{proof}
		
			\begin{proof}[Proof of Theorem~\ref{thm:algorithm}.]
				
				When the Lipschitz-correction step terminates, we have $\hat{L}_T\leq c \cdot L$ for each $T \in \mathcal{T}$. Therefore, $r_T \geq \nicefrac{\varepsilon}{2 (1+c)L}$. Consequently, the number of points needed in  $X_T$ to cover $T$ by balls of radius $r_T$ is bounded from above by a constant which depends on $\varepsilon$, $L$ and $c$ only. 
				Hence, any point added by the error-improvement step of Algorithm~\ref{alg:fitting} to $V$ has a distance of at least $r_{\min} = \nicefrac{\varepsilon}{2 (1+c)L}$ from the points of $V \cap T$. Since $\Omega$ is bounded, finitely many balls of radius $r_{\min}$ suffice to cover it completely, which proves finite convergence.
				
				When  Algorithm~\ref{alg:fitting} terminates, we can bound the maximum absolute difference of  $f$ and $\hat{f}$ as follows:
				\[
				\max_{x \in \Omega} |f(x) - \hat{f}(x)| = \max_{x \in T, T \in \mathcal{T}} |f(x) - \hat{f}(x)| \leq \max_{T \in \mathcal{T}} \varepsilon_T \overset{(1)}\leq \max_{T \in \mathcal{T}} \hat{\varepsilon}_T + r_T(L+\hat{L}_T) \overset{(2)}\leq \varepsilon/2 + \varepsilon/2 = \varepsilon,
				\]
				where (1) follows from Proposition~\ref{prop:sample-error} and (2) from the choice of $r_T$ and the terminating condition $\hat{\varepsilon}_{\max} =  \max_T \hat{\varepsilon}_T \leq \varepsilon/2$.
			\end{proof}

			\subsection{Piecewise-linear function fitting experiments for Section~\ref{sec:pwl_fit}}\label{app:fitting-exp}
			Further results for the PWL function fitting experiment are shown in Figures~\ref{fig:func1},~\ref{fig:func2} and~\ref{fig:func3}.
			To further support our method, we compared the number of triangles needed to approximate a set of bivariate test functions (see Table~\ref{tab:test-functions}) on the region $\left[0,1\right]^2$ within an error tolerance by our algorithm to the number of triangles needed by an equidistant grid triangulation, with grid diagonals chosen randomly.
			The error tolerances were $0.1,\,0.075,\,0.05,\,0.025$ and $0.01$.
			To determine the number of triangles needed for the grid triangulation, the number of grid points on the axes are increased until the measured error is below the tolerance.
			The results of the comparison are presented in Table~\ref{tab:trg-size}.
			
			\begin{table}[h!]
				\centering
				\caption{Test functions for PWL function fitting experiments.}\label{tab:test-functions}
				\begin{tabular}{l|l}
					\hline
					\textbf{Notation} & \textbf{Function} \\
					\hline
					$f_1$                 & $\exp\left(-5\left(\frac{\sqrt{(x-0.5)^2+(y-0.5)^2}}{1+0.3\sin\left(5\tan^{-1}\left(\frac{y-0.5}{x-0.5}\right)\right)}\right)^2\right)$ \\
					\hline
					$f_2$ & $\sin\left(6\pi x+ 0.5y\right)\exp\left(-10\left(\left(x-0.4\right)^2+\left(y-0.3\right)^2\right)\right)$\\ & $+ \cos\left(5\pi y+x\right)\exp\left(-12 \left(\left(x-0.7\right)^2 + \left(y-0.8\right)^2\right)\right)+ 0.1\sin\left(3\pi x y\right)$\\
					\hline
					$f_3$ & $\sin\left(3\pi x\right)\cos\left(\left(1-|y-0.5|\right)2\pi\right)\left(x+y\right)$ \\
					\hline
					$f_4$ & $\sin\left(50\sqrt{\left(y-0.5\right)^2+x^2}\right)\exp\left(-10\left(x^2+\left(y-0.5\right)^2\right)\right)$\\
					\hline
					$f_5$ & $\sin\left(5\pi x\right)\cos\left(5\pi y\right) + 0.5\sin\left(10\pi xy\right)+ 0.2\cos\left(15\left(x^2+y^2\right)\right)$\\
					\hline
				\end{tabular}
			\end{table}
			
			\begin{table}[h!]
				\small
				\caption{Number of simplices needed to approximate the test functions with maximal absolute error below tolerance levels by our algorithm and equidistant grid triangulations.}\label{tab:trg-size}
				\pgfplotstabletypeset[
				col sep=comma,
				header=true,
				string type,
				columns/Function/.style={column type=l, string type},
				columns/adaptive1/.style={column name=Adaptive, column type=r},
				columns/adaptive75/.style={column name=Adaptive, column type=r},
				columns/adaptive5/.style={column name=Adaptive, column type=r},
				columns/adaptive25/.style={column name=Adaptive, column type=r},
				columns/adaptive01/.style={column name=Adaptive, column type=r},
				columns/grid1/.style={column name=Grid, column type=r},
				columns/grid75/.style={column name=Grid, column type=r},
				columns/grid5/.style={column name=Grid, column type=r},
				columns/grid25/.style={column name=Grid, column type=r},
				columns/grid01/.style={column name=Grid, column type=r},
				every head row/.style={
					before row={
						\toprule
						Abs.err. & \multicolumn{2}{c}{0.1}
						& \multicolumn{2}{c}{0.075}
						& \multicolumn{2}{c}{0.05}
						& \multicolumn{2}{c}{0.025}
						& \multicolumn{2}{c}{0.01} \\
						\cmidrule(lr){2-3}
						\cmidrule(lr){4-5}
						\cmidrule(lr){6-7}
						\cmidrule(lr){8-9}
						\cmidrule(lr){10-11}
					},
					after row=\midrule
				},
				every last row/.style={
					after row=\bottomrule
				},
				row sep=\\,
				]{  Function,adaptive1,grid1,adaptive75,grid75,adaptive5,grid5,adaptive25,grid25,adaptive01,grid01\\
					$f_1$,102,338,138,512,204,1058,373,2592,1149,7442\\
					$f_2$,330,968,471,1352,690,2048,1544,4050,3860,10082\\
					$f_3$,254,578,324,800,502,1250,982,2450,2514,6050\\
					$f_4$,624,124002,920,152352,1408,168200,3081,215168,8634,>500000\\
					$f_5$,1172,7200,1453,9522,2401,14450,5092,28800,13169,72962\\
				}
			\end{table}

		\begin{figure}[h]

			\caption{Piecewise linear interpolation of $f_1$ on the $\left[0,1\right]^2$ domain with absolute error $\leq 0.005$.\label{fig:func1}}
		\end{figure}
		\begin{figure}[h]
			%
			\caption{Piecewise linear interpolation of $f_2$ on the $\left[0,1\right]^2$ domain with absolute error $\leq 0.025$.\label{fig:func2}}
		\end{figure}
		\begin{figure}[h]
			%
			\caption{Piecewise linear interpolation of $f_3$ on the $\left[0,1\right]^2$ domain with absolute error $\leq 0.025$.\label{fig:func3}}
		\end{figure}

		\section{Supplementary materials for Section~\ref{sec:mip_form}}
		\label{app:mip_form}
			\subsection{Proofs for Section~\ref{sec:rank_reduction}}
			\label{app:rank_reduction}
			
			\begin{proposition}\label{prop:lemma1}
				If $C$ is a minimal conflict set after splitting edge $(u,v)$ with new node $w$, $|C|\geq 3$, and $w \notin C$, then $C$ is a minimal conflict set before splitting edge $(u,v)$ as well.
			\end{proposition}
			\begin{proof}
				Let $C$ be as in the statement of the proposition. Suppose $C$ is not a conflict set before splitting edge $(u,v)$, that is, there exists a simplex $T$ such that $C\subseteq T$ before the split.
				Since $C$ is a conflict set after the split, this is only possible if $T$ is eliminated by the split. Hence $(u,v)$ must be an edge of $T$.
				There are two cases to consider. First suppose $\{u,v\} \subset C$. Then after  splitting  edge $(u,v)$, there is no simplex containing both $u$ and $v$, hence, $\{u,v\}$ is a two-element conflict set, and thus $C$ is not a minimal conflict set, a contradiction.
				Now suppose that at most one of $u$ and $v$ is in $C$. Since splitting edge $(u,v)$ of $T$ yields two simplices $T_1$ and $T_2$ such  that $u \in T_1$, $v \in T_2$, and $V(T_1) \setminus \{u\} = V(T_2)\setminus\{v\}$, $C$ is a subset of $T_1$ or $T_2$, a contradiction.
				
				Finally, we argue that $C$ was a minimal conflict set before the split. For suppose there exists $x\in C$ such that $C\setminus \{x\}$ is a conflict set before the split. We argue that $C\setminus\{x\}$ is a conflict set after the split as well, which contradicts that $C$ is a minimal conflict set after the split.
				So suppose $C\setminus\{x\}$ is not a conflict set after splitting edge $(u,v)$. Then there exists a simplex $T$ containing $C\setminus \{x\}$. 
				If $T$ is not obtained by splitting edge $(u,v)$, then $C\setminus\{x\}$ is not a conflict set before the split, a contradiction.
				If $T$ is obtained by splitting edge $(u,v)$, then there exists a simplex $T'$ containing edge $(u,v)$ such that $T \subset T'$. Hence, $C\setminus  \{x\} \subseteq T \subset T'$. This contradicts our assumption that  $C\setminus\{x\}$ is a conflict set before the split. Therefore, $C$ was a minimal conflict set before the split as well.
			\end{proof}
			
			\begin{proof}[Proof of Proposition~\ref{prop:splitting2}]
				Proposition~\ref{prop:lemma1} implies that if $C$ is a newly formed conflict set after the split of the edge $(u,v)$ by node $w$, then $w\in C$.
				Since $C$ is a minimal conflict set after the split, $C\setminus\{w\}$ cannot be a conflict set. Clearly, if $\{u,v\} \subset C$, then $C$ cannot be a minimal conflict set, since $\{u,v\}$ is a minimal conflict set after the split.
				If $C\setminus\{w\} \cup \{u\}$ 
				is a minimal conflict set, then $u\notin C$ and $C\setminus\{w\} \cup \{u\}$ is a minimal conflict set before the split by the previous lemma.
				By a similar argument, if $C\setminus\{w\} \cup \{v\}$ is a minimal conflict set after the split, then $v\notin C$ and $C\setminus\{w\} \cup \{v\}$ is a minimal conflict set before the split. 
				
				So assume  neither of the cases \ref{split:case1} or \ref{split:case2} hold. 
				First, we argue that $\{u,v\} \cap C = \emptyset$.
				If $u \in C$, then after the split, the minimal conflict set $\{u,v\}$ is a subset of $C\setminus\{w\} \cup \{v\}$, hence, $C\setminus\{w\} \cup \{v\}$ is a conflict set, a contradiction.
				A similar argument shows that $v \notin C$.

				Now we probe that case \ref{split:case3} holds. For the sake of a contradiction,
				suppose $T$ is a simplex such that $C\setminus\{w\} \cup \{u,v\} \subseteq T$ before the split. Since  $\{u,v\}\subset T$,  $T$ is split into two simplices $T^s_u$ and $T^s_v$ when splitting edge $(u,v)$ by $w$ such that $u \in T^s_u$ and $v \in T^s_v$. Therefore, $C\setminus\{w\} \cup \{u\} \subset T^s_u$ and  $C\setminus\{w\} \cup \{v\} \subset T^s_v$.
				Then, $C \subset C\cup \{u\} \subset T^s_u$, since $w \in T^s_u$. This contradicts the assumption that $C$ is a conflict set after splitting edge $(u,v)$ by $w$.
				Hence, $C\setminus\{w\} \cup \{u,v\}$ is a conflict set before the split.
				It is a minimal conflict set by the above assumption, and case \ref{split:case3} follows.
			\end{proof}
			
			\begin{proof}[Proof of Proposition~\ref{prop:splitting3}]
				Consider any $w' \in C\setminus\{w\}$.
				Since $C$ is a minimal conflict set after the split and $|C| \geq 3$, there exists a simplex $T$  containing both $w$ and $w'$.
				Since $w \in C$, $T$ was created by  splitting edge $(u,v)$ by $w$.
				Then there existed a simplex $T'$ before the split containing both $u$ and $v$ such that $V(T)\setminus \{w\} \subset V(T')$.
				Then $w', u, v \in T'$. Let $T^s_u$ and $T^s_v$ be the simplices obtained from $T'$ by splitting edge $(u,v)$ by $w$ such that $u \in T^s_u$, and $v \in T^s_v$. Then $u,w'\in T^s_u$ and $v,w'\in T^s_v$ as claimed.
			\end{proof}

		\subsection{Proofs for Section~\ref{sec:high-rank-conflict}}\label{app:high-rank-conflict}
		\begin{proof}[Proof of Proposition~\ref{prop:blocking-hg}]
			To prove~\ref{prop:blocking-hg:a}, consider any hyperedge $C=\left\{v_1,\ldots,v_k\right\}$ of $\Hc_{\mathcal{S}}$ with $k\geq 3$.
			Since $C$ is a minimal infeasible set, there exist sets $S_i, S_j \in \mathcal{S}$ such that  $C\setminus\left\{v_i\right\}\subseteq S_i$ and $C\setminus\left\{v_j\right\}\subseteq S_j$ and $C  \subseteq S_i\cup S_j$ for any pair of distinct indices $i$ and $j$.
			Then $\left\{S_i,S_j\right\}$ is a minimal blocking set, hence, a hyperedge of the blocking hypergraph.
			
			As for property~\ref{prop:blocking-hg:b}, to cover a $k$-element subset of $V$, one needs at most $k$ sets $S_i$ from $\mathcal{S}$.
		\end{proof}

		\begin{proof}[Proof of Proposition~\ref{prop:cdc-color}]
			By (\ref{eq:colormip}), there  exists a unique color $c$ such that $z_c = 1$, and $z_{c'} = 0$ for all $c' \neq c$. Hence, $\supp{\lambda}$ is a subset of the vertices of the polytopes in $K_c$. However, $K_c$ does not contain a blocking set of polytopes by Proposition~\ref{prop:no_blocking}. Therefore, $\supp{\lambda}$ does not contain a minimal infeasible set of size at least 3, and the statement follows.
		\end{proof}

		\subsection{Coloring the blocking hypergraph for Section~\ref{sec:high-rank-conflict}}\label{app:coloring}
		Below we present a method for finding a coloring of the nodes of the blocking hypergraph $\Hb_{\mathcal{S}} = (\mathcal{S},\mathcal{E}^b_\mathcal{S})$ as defined in Section~\ref{sec:mip_form}.
		Let $\mathcal{E}^{2}_\mathcal{S}=\{B\in \mathcal{E}^b_\mathcal{S}\ :\ |B|=2\}$ 
		be the subset of hyperedges of rank $2$.
		Observe, that $\mathcal{E}^2_\mathcal{S}\neq\emptyset$ by Proposition~\ref{prop:blocking-hg}~\ref{prop:blocking-hg:a}.
		
		First, the problem is reduced to the coloring of $\mathcal{H}^2_\mathcal{S}=(\mathcal{S},\mathcal{E}^2_\mathcal{S})$, which is solved greedily as an ordinary graph coloring problem.
		Suppose a coloring with $q$ colors is found. 
		Then the  hypergraph colorability problem with $q$ colors is modeled as a boolean satisfiability problem (SAT), by a straightforward formulation of the hypergraph coloring problem.
		Let $X_{S,c}$ denote the variable corresponding to coloring simplex $S\in\mathcal{S}$ with color $c\in\range{q}$.
		Then every assignment of truth values of $X$ that satisfies the formula
		\begin{equation}
			\left(\bigwedge_{S\in\mathcal{S}}\bigvee_{c\in\range{q}}X_{S,c}\right)\land\left(\bigwedge_{S\in\mathcal{S}}\bigwedge_{c\neq c'\in\range{q}}\left(\lnot X_{S,c}\lor\lnot X_{S,c'}\right)\right)\land\left(\bigwedge_{B\in\mathcal{E}^b_\mathcal{S}}\bigwedge_{c\in\range{q}}\bigvee_{S\in B}\lnot X_{S,c}\right)
		\end{equation}
		define a coloring of the blocking hypergraph $\Hb_{\mathcal{S}}$.
		The first two terms encode that each polyherdon is colored with exactly one color, while the last term encodes that each hyperedge has at least one polyhedron that is colored with at least two colors.
		We have $\rank{\Hb_\mathcal{S}}\leq d+1$ by Theorem~\ref{thm:conflict-hg} and Proposition~\ref{prop:blocking-hg}~\ref{prop:blocking-hg:b}, hence, computing $\Hb_\mathcal{S}$ can be done by checking all subsets of $\mathcal{S}$ of cardinality at most $d+1$.
		The formula is solved using specialized software, see~\citep{audemard_glucose_2018}. 
		If a feasible coloring with $q$ colors is found, then $q$ is decreased by $1$, and the procedure is repeated until the formula becomes infeasible.
		Otherwise, $q$ is increased by $1$, and the procedure is repeated  until a feasible solution is found.

		\begin{example}
			Suppose that $d=2$ and $\mathcal{S} = \left(S_1,S_2,S_3,S_4\right)$ is the simplicial partition depicted in Figure~\ref{fig:example:a}. Then $\Hc_{\mathcal{S}} = (V,\mathcal{E}_\mathcal{S})$ with $V=\left\{u,v_1,v_2,v_3,w\right\}$ and $\mathcal{E}_\mathcal{S} = \left\{\{u,w\},\,\{v_3,w\},\,\{v_1,v_2,v_3\}\right\}$.
			The conflict graph of $\mathcal{S}$ is $G^c_\mathcal{S}=(V,\{\{u,w\},\,\{v_3,w\}\})$, and $A_1 = \{w\}$ and $B_1=\{u,v_3\}$ induces a biclique  of $G^c_\mathcal{S}$.
			Then hyperedge $\{v_1,v_2,v_3\}$ is not excluded from the support of $\lambda$ by system~(\ref{eq:pairwise-IB}).
			The blocking hypegraph of $\mathcal{S}$ is $\Hb_{\mathcal{S}}=(\mathcal{S}, \mathcal{E}^b_{\mathcal{S}})$, with	$\mathcal{E}^b_{\mathcal{S}} = \left\{\{S_1,S_2\}, \{S_1,S_3\}, \{S_2,S_3\}, \{S_2,S_4\},\{S_3,S_4\}\right\}$ is depicted in Figure~\ref{fig:example:b}.
			Observe that the rank of the blocking hypergraph is 2, hence, it can be colored as an ordinary graph.
			There exists a proper coloring of $\Hb_{\mathcal{S}}$ with 3 colors as follows.
			Let $K_1 = \left\{S_1,S_4\right\},\, K_2 = \left\{S_2\right\}$ and $K_3 = \left\{S_3\right\}$, see Figure~\ref{fig:example:c}.
			The coloring patterns of the nodes in $V$ are $\pi_u = \left\{1,2,3\right\},\, \pi_{v_1} = \left\{1,2\right\},\, \pi_{v_2} = \left\{1,3\right\},\, \pi_{v_3} = \left\{2,3\right\},\, \pi_w = \left\{1\right\}$.
			Then the MIP formulation of (\ref{eq:disj_constr}) over $\mathcal{S}$ is
			\begin{align*}
				&\lambda_w \leq y_1,\quad
				\lambda_u+\lambda_{v_3} \leq 1-y_1,\\
				&\lambda_u \leq z_1+z_2+z_3, \quad
				\lambda_{v_1} \leq  z_1+z_2, \quad
				\lambda_{v_2} \leq  z_1+z_3,\quad
				\lambda_{v_3} \leq  z_2+z_3,\quad
				\lambda_w \leq  z_1,\quad z_1+z_2+z_3 = 1,\\
				&\lambda \in \Delta^V,\quad
				y_1,z_1,z_2,z_3\in\left\{0,1\right\}.
			\end{align*}
		\end{example}
		
		\begin{figure}[h]
			\small
			\centering
			\subcaptionbox{\label{fig:example:a}}
			{
				\begin{tikzpicture}[scale=0.9]
					\coordinate (v1) at (1,3);
					\node at (v1) [above = 1mm of v1] {$v_1$};
					\coordinate (v2) at (3,1);
					\node at (v2) [below = 1mm of v2] {$v_2$};
					\coordinate (v3) at (0,0);
					\node at (v3) [left = 1mm of v3] {$v_3$};
					\coordinate (u) at (1.3,1.3);
					\node at (u) [below = 1mm of u] {$u$};
					\coordinate (w) at (3.3,3.3);
					\node at (w) [above = 1mm of w] {$w$};
					\draw (v1) -- (v2) -- (v3) -- (v1);
					\draw (v1) -- (u);
					\draw (v2) -- (u);
					\draw (v3) -- (u);
					\draw (v1) -- (w) -- (v2);
					\draw[red, opacity=.5, fill=red!50, rounded corners=5pt] ($(v2)+(.3,0)$) -- ($(v1)+(0,.3)$) -- ($(v3)+(-.3,0)$) -- ($(v3)+(0,-.3)$) -- ($(v1)+(0.1,-.3)$)  -- ($(v2)+(-.3,0)$)  -- ($(v2)+(0,-.3)$) --cycle;
					\path[dashed, red] (u)  edge   [bend right=15]   node {} (w);
					\path[dashed, red] (v3)  edge   [bend left=15]   node {} (w);
				\end{tikzpicture}
			}
			\hfill
			\subcaptionbox{\label{fig:example:b}}
			{
				\small
				\centering
				\begin{tikzpicture}[scale=0.9]
					\coordinate (v1) at (1,3);
					\coordinate (v2) at (3,1);
					\coordinate (v3) at (0,0);
					\coordinate (u) at (1.3,1.3);
					\coordinate (w) at (3.3,3.3);
					\node at (v1) [above = 1mm of v1] {$v_1$};
					\node at (v2) [below = 1mm of v2] {$v_2$};
					\node at (v3) [left = 1mm of v3] {$v_3$};
					\node at (w) [above = 1mm of w] {$w$};
					\draw (v1) -- (v2) -- (v3) -- (v1);
					\draw (v1) -- (u);
					\draw (v2) -- (u);
					\draw (v3) -- (u);
					\draw (v1) -- (w) -- (v2);
					\node (T1) at ($.33*(v1)+.33*(v2)+.33*(u)$) {$S_1$};
					\node (T2) at ($.33*(v1)+.33*(v3)+.33*(u)$) {$S_2$};
					\node (T3) at ($.33*(v3)+.33*(v2)+.33*(u)$) {$S_3$};
					\node (T4) at ($.33*(v1)+.33*(v2)+.33*(w)$) {$S_4$};
					\path[green!60!black] (T2)  edge   [bend left=20]   node {} (T4);
					\path[green!60!black] (T3)  edge   [bend right=20]   node {} (T4);
					\path[green!60!black] (T1)  edge   [bend right=20]   node {} (T2);
					\path[green!60!black] (T1)  edge   [bend right=20]   node {} (T3);
					\path[green!60!black] (T2)  edge   [bend right=20]   node {} (T3);
				\end{tikzpicture}
			}
			\hfill
			\subcaptionbox{\label{fig:example:c}}{
				\small
				\centering
				\begin{tikzpicture}[scale=0.9]
					\coordinate (v1) at (1,3);
					\node at (v1) [above = 1mm of v1] {$v_1$};
					\coordinate (v2) at (3,1);
					\node at (v2) [below = 1mm of v2] {$v_2$};
					\coordinate (v3) at (0,0);
					\node at (v3) [left = 1mm of v3] {$v_3$};
					\coordinate (u) at (1.3,1.3);
					\coordinate (w) at (3.3,3.3);
					\node at (w) [above = 1mm of w] {$w$};
					\draw[opacity=.5, fill=blue!70!black] (u) -- (v1) -- (w)--(v2);
					\draw[opacity=.5, fill=red!70!black] (u) -- (v1) -- (v3);
					\draw[opacity=.5, fill=yellow!70!black] (u) -- (v3) -- (v2);
					\draw (v1) -- (v2) -- (v3) -- (v1);
					\draw (v1) -- (u);
					\draw (v2) -- (u);
					\draw (v3) -- (u);
					\draw (v1) -- (w) -- (v2);
					\node (T1) at ($.33*(v1)+.33*(v2)+.33*(u)$) {$S_1$};
					\node (T2) at ($.33*(v1)+.33*(v3)+.33*(u)$) {$S_2$};
					\node (T3) at ($.33*(v3)+.33*(v2)+.33*(u)$) {$S_3$};
					\node (T4) at ($.33*(v1)+.33*(v2)+.33*(w)$) {$S_4$};
					\node at (u) [below = .5mm of u] {$u$};
				\end{tikzpicture}
			}
		\caption{Example of a polyhedral partition with conflict hypergraph of rank 3.\label{fig:example}}
		\end{figure}
		
		\subsection{Maximum weight biclique heuristic for planar graphs with non-negative edge-weights}\label{app:biclique_heur}
		We describe a randomized algorithm for finding bicliques with large weight if $G$ is planar.
		It is  assumed that a planar embedding of $G$ is known.
		This setting arises naturally when $G$ equals the complement of the conflict graph $G^c_\mathcal{S}$ for a set-system $\mathcal{S}$.
		
		The heuristic is based on the following observation about the structure of bicliques in a planar graph.
		A subset $S$ of nodes of a connected graph $G=(V,E)$ is a {\em vertex cut}, if removing the subset of nodes $S$ from $V$ disconnects the graph into at least two connected components.
		
		Suppose $S$ is a vertex cut in $G$,  graph $G'=(V\setminus S,E')$ with $E'$ being the subset of edges of $G$ spanned by $V\setminus S$ has $k\geq 2$ connected components, and the node sets of those components are $V_1,\ldots,V_k$.
		It's easy to see that  for any $i,j\in\range{k},\,i\neq j$ and $E_{ij}=\left\{\{u,v\}\vert u\in V_i,\, v\in V_j\right\}$,  $(V_i\cup V_j,E_{ij})$ is a biclique in $\overline{G}$.

		Let $G^*$ denote the dual graph of $G$, and $G^{*}_-$ is obtained from $G^*$ by removing the node corresponding to the unbounded face in the planar embedding of $G$.
		Graph $G^{*}_-$ is also planar, and a planar embedding of it can be derived from that of $G$.
		Also, by F\'ary's theorem \citep{fary1948straight}, $G^{*}_-$ can be drawn such that its edges are straight line segments.
		Let $U\subseteq V$ denote those nodes of $G$  incident to the unbounded face of its embedding in the plane.
		\begin{proposition}\label{prop:path}
		Suppose $\ell\subset\R^2$ is a line in the plane that intersects the planar embedding of $G^{*}_-$ in a positive finite number of points and the points of intersection are not vertices of $G^{*}_-$.
		Let $E^*_\ell$ denote the edges of $G^{*}_-$ intersected by $\ell$ and $E_\ell$ denote the corresponding edges of $G$.
		Then $E_\ell$ forms a path in $G$ with the two endpoints in $U$.
		\begin{proof}{Proof.}
			Suppose that $\ell$ is directed, and let $e^*_1$ denote the first edge of $G^*_-$ intersected by $\ell$, and $P^*_1$ the first face of $G^*_-$ crossed by $\ell$.
			Then $e^*_1$ intersects a unique edge $e_1$ of $G$ in the given planar embedding.
			Let $v_1$ denote the node of $G$ corresponding to $P^*_1$ and $v_0$  the other endpoint of $e_1$.
			Then $v_1$ lies in the interior of face $P^*_1$ (in the embedding), $v_0$ lies outside it, and there is no face of $G^*_-$ that contains $v_0$ in its interior.
			Hence, $v_0\in U$.
			Let $P^*_i$ and $P^*_{i+1}$ denote the $i$th and $(i+1)$th face of $G^*_-$ crossed by $\ell$, and $v_i,v_{i+1}$ the corresponding nodes of $G$.
			Then $v_iv_{i+1}$ is an edge of $G$, corresponding to the boundary edge of $P^*_i$ and $P^*_{i+1}$.
			Let $P^*_k$ denote the last face and $e^*_k$ the last edge of $G^*_-$ crossed by $\ell$.
			Similarly to $e^*_1$, the dual of $e^*_k$ has and endpoint (denoted by $v_k$) in $P^*_k$, and the other endpoint $v_{k+1}\in U$.
			Then the nodes $v_0,v_1,\ldots,v_k,v_{k+1}$ along with the edges between consecutive points form a path in $G$ with $v_0,v_{k+1}\in U$.
		\end{proof}
		\end{proposition}

		Based on the above observations, a randomized heuristic algorithm is proposed.
		The heuristic is illustrated in Figure~\ref{fig:heur}.
		Figures~\ref{fig:heur:a} and~\ref{fig:heur:b} show the complement $G$ of a conflict graph $G^c_\mathcal{S}$ and its dual $G^*_-$ and their planar embedding.
		
		First $n$ random straight lines ($\ell_1,\ldots,\ell_n$) that intersect the planar embedding of $G^*_-$ are generated.
		Then for each line $\ell$ (dashed blue line in Figure~\ref{fig:heur:c}), the set of edges of $G^*_-$ crossed by $\ell$ are computed (thick red edges in Figure~\ref{fig:heur:d}), as well as the edges of $G$ corresponding to these edge (solid blue edges in Figure~\ref{fig:heur:e}).
		By Proposition~\ref{prop:path}, the selected edges of $G$ constitute  a simple path in $G$ and consequently a vertex cut of $G$, which in turn corresponds to a biclique in $\overline{G}$ (the solid and empty green colored nodes and curved edges in Figure~\ref{fig:heur:f}, the empty node constitutes one part of the biclique and the solid nodes the other).
		After computing the biclique for each $\ell_i$, the one with the maximum weight is chosen. Detailed computation results for the above techniques are presented in Section~\ref{app:bicliqueCover:exp}.
		
		The above method can be further refined to identify more complex bicliques by iteratively cutting $G^*_-$ and its resulting subgraphs, then merging the connected components into just two parts, which are subsequently tested for maximal weight bicliques.
		\begin{figure}[ht]
		\subcaptionbox{\label{fig:heur:a}}{
			\centering
			\begin{tikzpicture}[scale=2.5]
				\coordinate (a) at (0,0);
				\coordinate (b) at (0,1);
				\coordinate (c) at (1,0);
				\coordinate (d) at (1,1);
				\coordinate (e) at (0.3,0.2);
				\coordinate (f) at (0.4,0.7);
				\coordinate (g) at (.6,0);
				\draw (a) -- (b) -- (d) -- (c) -- (a);
				\draw (a) -- (f);
				\draw (a) -- (e);
				\draw (b) -- (f);
				\draw (c) -- (e);
				\draw (c) -- (f);
				\draw (d) -- (f);
				\draw (g) -- (e);
				\draw (e) -- (f);
			\end{tikzpicture}
		}
		\hfill
		\subcaptionbox{\label{fig:heur:b}}{
			\centering
			\begin{tikzpicture}[scale=2.5]
				\coordinate (a) at (0,0);
				\coordinate (b) at (0,1);
				\coordinate (c) at (1,0);
				\coordinate (d) at (1,1);
				\coordinate (e) at (0.3,0.2);
				\coordinate (f) at (0.4,0.7);
				\coordinate (g) at (.6,0);
				\draw (a) -- (b) -- (d) -- (c) -- (a);
				\draw (a) -- (f);
				\draw (a) -- (e);
				\draw (b) -- (f);
				\draw (c) -- (e);
				\draw (c) -- (f);
				\draw (d) -- (f);
				\draw (g) -- (e);
				\draw (e) -- (f);
				\coordinate (abf) at ($.33*(a)+.33*(b)+.33*(f)$);
				\coordinate (aef) at ($.33*(a)+.33*(e)+.33*(f)$);
				\coordinate (aeg) at ($.33*(a)+.33*(e)+.33*(g)$);
				\coordinate (bfd) at ($.33*(d)+.33*(b)+.33*(f)$);
				\coordinate (cfd) at ($.33*(d)+.33*(c)+.33*(f)$);
				\coordinate (efc) at ($.33*(e)+.33*(c)+.33*(f)$);
				\coordinate (gec) at ($.33*(g)+.33*(e)+.33*(c)$);
				\draw[red] (abf) -- (aef) -- (aeg) -- (gec) -- (efc) -- (cfd) -- (bfd) -- (abf);
				\draw[red] (efc) -- (aef);
			\end{tikzpicture}
		}
		\hfill
		\subcaptionbox{\label{fig:heur:c}}{
			\centering
			\begin{tikzpicture}[scale=2.5]
				\coordinate (a) at (0,0);
				\coordinate (b) at (0,1);
				\coordinate (c) at (1,0);
				\coordinate (d) at (1,1);
				\coordinate (e) at (0.3,0.2);
				\coordinate (f) at (0.4,0.7);
				\coordinate (g) at (.6,0);
				\draw (a) -- (b) -- (d) -- (c) -- (a);
				\draw (a) -- (f);
				\draw (a) -- (e);
				\draw (b) -- (f);
				\draw (c) -- (e);
				\draw (c) -- (f);
				\draw (d) -- (f);
				\draw (g) -- (e);
				\draw (e) -- (f);
				\coordinate (abf) at ($.33*(a)+.33*(b)+.33*(f)$);
				\coordinate (aef) at ($.33*(a)+.33*(e)+.33*(f)$);
				\coordinate (aeg) at ($.33*(a)+.33*(e)+.33*(g)$);
				\coordinate (bfd) at ($.33*(d)+.33*(b)+.33*(f)$);
				\coordinate (cfd) at ($.33*(d)+.33*(c)+.33*(f)$);
				\coordinate (efc) at ($.33*(e)+.33*(c)+.33*(f)$);
				\coordinate (gec) at ($.33*(g)+.33*(e)+.33*(c)$);
				\draw[red] (abf) -- (aef) -- (aeg) -- (gec) -- (efc) -- (cfd) -- (bfd) -- (abf);
				\draw[red] (efc) -- (aef);
				\coordinate (l1) at (.2,1);
				\coordinate (l2) at (.4,0);
				\draw[blue, dashed] (l1) -- (l2);
			\end{tikzpicture}
		}
		\subcaptionbox{\label{fig:heur:d}}{
			\centering
			\begin{tikzpicture}[scale=2.5]
				\coordinate (a) at (0,0);
				\coordinate (b) at (0,1);
				\coordinate (c) at (1,0);
				\coordinate (d) at (1,1);
				\coordinate (e) at (0.3,0.2);
				\coordinate (f) at (0.4,0.7);
				\coordinate (g) at (.6,0);
				\draw (a) -- (b) -- (d) -- (c) -- (a);
				\draw (a) -- (f);
				\draw (a) -- (e);
				\draw (b) -- (f);
				\draw (c) -- (e);
				\draw (c) -- (f);
				\draw (d) -- (f);
				\draw (g) -- (e);
				\draw (e) -- (f);
				\coordinate (abf) at ($.33*(a)+.33*(b)+.33*(f)$);
				\coordinate (aef) at ($.33*(a)+.33*(e)+.33*(f)$);
				\coordinate (aeg) at ($.33*(a)+.33*(e)+.33*(g)$);
				\coordinate (bfd) at ($.33*(d)+.33*(b)+.33*(f)$);
				\coordinate (cfd) at ($.33*(d)+.33*(c)+.33*(f)$);
				\coordinate (efc) at ($.33*(e)+.33*(c)+.33*(f)$);
				\coordinate (gec) at ($.33*(g)+.33*(e)+.33*(c)$);
				\draw[red] (abf) -- (aef) -- (aeg) -- (gec) -- (efc) -- (cfd) -- (bfd) -- (abf);
				\draw[red, line width=0.5mm,] (abf) -- (bfd);
				\draw[red, line width=0.5mm,] (efc) -- (aef);
				\draw[red, line width=0.5mm,] (aeg) -- (gec);
				\coordinate (l1) at (.2,1);
				\coordinate (l2) at (.4,0);
				\draw[blue, dashed] (l1) -- (l2);
			\end{tikzpicture}
		}
		\hfill
		\subcaptionbox{\label{fig:heur:e}}{
			\centering
			\begin{tikzpicture}[scale=2.5]
				\coordinate (a) at (0,0);
				\coordinate (b) at (0,1);
				\coordinate (c) at (1,0);
				\coordinate (d) at (1,1);
				\coordinate (e) at (0.3,0.2);
				\coordinate (f) at (0.4,0.7);
				\coordinate (g) at (.6,0);
				\draw (a) -- (b) -- (d) -- (c) -- (a);
				\draw (a) -- (f);
				\draw (a) -- (e);
				\draw (b) -- (f);
				\draw (c) -- (e);
				\draw (c) -- (f);
				\draw (d) -- (f);
				\draw (g) -- (e);
				\draw (e) -- (f);
				\coordinate (abf) at ($.33*(a)+.33*(b)+.33*(f)$);
				\coordinate (aef) at ($.33*(a)+.33*(e)+.33*(f)$);
				\coordinate (aeg) at ($.33*(a)+.33*(e)+.33*(g)$);
				\coordinate (bfd) at ($.33*(d)+.33*(b)+.33*(f)$);
				\coordinate (cfd) at ($.33*(d)+.33*(c)+.33*(f)$);
				\coordinate (efc) at ($.33*(e)+.33*(c)+.33*(f)$);
				\coordinate (gec) at ($.33*(g)+.33*(e)+.33*(c)$);
				\draw[red] (abf) -- (aef) -- (aeg) -- (gec) -- (efc) -- (cfd) -- (bfd) -- (abf);
				\draw[red] (efc) -- (aef);
				\coordinate (l1) at (.2,1);
				\coordinate (l2) at (.4,0);
				\draw[blue, dashed] (l1) -- (l2);
				\draw[blue, thick] (b) -- (f) -- (e) -- (g);
			\end{tikzpicture}
		}
		\hfill
		\subcaptionbox{\label{fig:heur:f}}{
			\centering
			\begin{tikzpicture}[scale=2.5]
				\coordinate (a) at (0,0);
				\coordinate (b) at (0,1);
				\coordinate (c) at (1,0);
				\coordinate (d) at (1,1);
				\coordinate (e) at (0.3,0.2);
				\coordinate (f) at (0.4,0.7);
				\coordinate (g) at (.6,0);
				\draw (a) -- (b) -- (d) -- (c) -- (a);
				\draw (a) -- (f);
				\draw (a) -- (e);
				\draw (b) -- (f);
				\draw (c) -- (e);
				\draw (c) -- (f);
				\draw (d) -- (f);
				\draw (g) -- (e);
				\draw (e) -- (f);
				\coordinate (abf) at ($.33*(a)+.33*(b)+.33*(f)$);
				\coordinate (aef) at ($.33*(a)+.33*(e)+.33*(f)$);
				\coordinate (aeg) at ($.33*(a)+.33*(e)+.33*(g)$);
				\coordinate (bfd) at ($.33*(d)+.33*(b)+.33*(f)$);
				\coordinate (cfd) at ($.33*(d)+.33*(c)+.33*(f)$);
				\coordinate (efc) at ($.33*(e)+.33*(c)+.33*(f)$);
				\coordinate (gec) at ($.33*(g)+.33*(e)+.33*(c)$);
				\coordinate (l1) at (.2,1);
				\coordinate (l2) at (.4,0);
				\node[circle,draw=black!60!green, thick, fill=white, inner sep=2pt,minimum size=3pt] (u) at (0,.03) {};
				\node[circle,draw=black!60!green, thick, fill=black!60!green, inner sep=2pt,minimum size=3pt] (v1) at (1,.97) {};
				\node[circle,draw=black!60!green, thick, fill=black!60!green, inner sep=2pt,minimum size=3pt] (v2) at (1,.03) {};
				\path[black!60!green, thick] (u)  edge   [bend left=25]   node {} (v1);
				\path[black!60!green, thick] (u)  edge   [bend left=30]   node {} (v2);
			\end{tikzpicture}
		}
		\caption{Randomized heuristic for maximum weight biclique in a co-planar graph.\label{fig:heur}}
		\end{figure}

		\subsection{Computational experiments for Section~\ref{sec:pairwise-conflicts}}\label{app:bicliqueCover:exp}
		In this section we describe computational results for two variants of Algorithm~\ref{alg:biclique_cover}.
		In the first one, termed 'Pure IP', we find maximum weight bicliques by solving MIP (\ref{eq:maxBiclique}) in each iteration.
		In the second variant, termed 'Geom+IP', we applied the randomized geometric heuristic for four iterations and then went on with solving IP~\eqref{eq:maxBiclique}.
		
		We compared these methods on random triangulations of a rectangular region in the plane of three different sizes, 10 instances each, described in Table~\ref{tab:randomSize}.
		Columns $|E(G^c_\mathcal{S})|$ and $\lceil \log_2|E(G^c_\mathcal{S})|\rceil$ contain the number of edges of the conflict graph and its base 2 logarithm.
		The random triangulations were generated in three size classes (small, medium and large), with each size class containing 10 instances.
		The size of the instances are summarized in Table~\ref{tab:randomSize}.
		The triangulations were generated by taking an initial rectangle and its four corner points, and adding a pre-specified number of points generated uniformly at random within the rectangle, then taking a Delaunay-triangulation of the points.
		
		A time limit of \SI{100}{\second} were set for solving IP~\eqref{eq:maxBiclique}.
		For the randomized heuristic, in each of the first four iterations 1000 random lines were generated.
		The left hand side of Figure~\ref{fig:avgBiclique} shows the average number of bicliques found for the three problem sizes by the strategies.
		The horizontal lines signify $\log_2|E(G^c_\mathcal{S})|$ for the different triangle sizes.
		Observe that many of the times a biclique cover of size $\log_2|E(G^c_\mathcal{S})|$ can be found.
		In terms of runtime (on the  right pane of Figure~\ref{fig:avgBiclique}), strategy Geom+IP could beat strategy Pure IP on the large instances, lost to it on the small instances, and produced similar runtimes on the medium size instances.
		The reason for that is solving and proving optimality for IP~(\ref{eq:maxBiclique}) is harder when there is a large number of edges with non-zero weights, and strategy Geom+IP used the fast, randomized heuristic in this case.
		The experiments show that  our greedy heuristic performs well on co-planar graphs.
		
		\begin{table}
		\centering
		\small
		\caption{Number of triangles, ground points and size of biclique cover for different sizes of random triangulations.
			\label{tab:randomSize}}
		\footnotesize
		\pgfplotstabletypeset[
		col sep=comma,
		columns/Type/.style={string type}, 
		columns/Size/.style={string type}, 
		columns={ Size, Ntrg, Npts, Nedge, LogEdge}, 
		every head row/.style={before row=\toprule, after row=\midrule}, 
		every last row/.style={after row=\bottomrule}, 
		columns/Type/.style={string type, column type=l}, 
		columns/Size/.style={string type, column type=r}, 
		columns/Ntrg/.style={column type=r, column name=\#triangles}, 
		columns/Npts/.style={column type=r, column name=\#points}, 
		columns/Nbiclique/.style={string type, column type=r, column name=Size of biclique cover}, 
		columns/Ncolor/.style={string type, column type=r, column name=Number of colors}, 
		columns/Nedge/.style={column type=r, column name=$|E(G^c_\mathcal{S})|$}, 
		columns/LogEdge/.style={column type=r, column name=$\lceil \log_2|E(G^c_\mathcal{S})|\rceil$}, 
		every head row/.style={
			before row={
				\toprule
			},
			after row=\midrule
		},
		row sep=\\,
		]{
			Type, Size, Ntrg, Npts, Nbiclique, Ncolor, Nedge, LogEdge\\
			random, small, 26, 18, {7-8}, {1-3}, 110, 7\\
			, medium, 64, 41, {9-10}, {1-3}, 716, 10\\
			, large, 126, 78, {11-12}, {1-6}, 2800, 12\\
		}
		\end{table}

		\begin{figure}[h]
			\begin{tikzpicture}[scale=.8]
				\begin{groupplot}[
					group style={
						group size=2 by 1,
						x descriptions at=edge bottom,
						horizontal sep=50pt,
					},
					]
					\nextgroupplot[
					ybar, 
					bar width=22pt, 
					symbolic x coords={Strategy\#1,Strategy\#4}, 
					xticklabels={Geom+IP, Pure IP},
					xtick style={draw=none},
					xtick=data, 
					xlabel={}, ylabel={Avg. \#bicliques},
					width=.5\textwidth,
					height=6cm,
					grid=major,
					ymin=5, 
					ymax=16,
					ytick={5,7,10,12,15},
					nodes near coords, 
					every node near coord/.append style={anchor=north, yshift=30pt}, 
					enlarge x limits=.5, 
					]
					\addplot+[error bars/.cd, y dir=both, y explicit, error mark=triangle*] 
					table[row sep=\\,  y error plus=Max, y error minus=Min] {
						Category   Mean Min Max \\
						Strategy\#1  7.8  .8 1.2 \\
						Strategy\#4  7.5  .5 .5 \\
					};
					\addlegendentry{Small} 
					\addplot+[error bars/.cd, y dir=both, y explicit, error mark=triangle*] 
					table[row sep=\\,  y error plus=Max, y error minus=Min] {
						Category   Mean Min Max \\
						Strategy\#1  10   1    1 \\
						Strategy\#4  9.7   .7   .3 \\
					};
					\addlegendentry{Medium} 
					\addplot+[error bars/.cd, y dir=both, y explicit, error mark=triangle*] 
					table[row sep=\\,  y error plus=Max, y error minus=Min] {
						Category   Mean Min Max \\
						Strategy\#1  11.4  .4  .6\\
						Strategy\#4  11.4  .4  .6\\
					};
					\addlegendentry{Large} 
					\coordinate (small) at (axis cs:Strategy\#1,7);
					\coordinate (medium) at (axis cs:Strategy\#1,10);
					\coordinate (large) at (axis cs:Strategy\#1,12);
					\coordinate (O1) at (rel axis cs:0,0);
					\coordinate (O2) at (rel axis cs:1,0);
					\draw [blue,sharp plot,dashed] (small -| O1) -- (small -| O2);
					\draw [red,sharp plot,dashed] (medium -| O1) -- (medium -| O2);
					\draw [brown!60!black,sharp plot,dashed] (large -| O1) -- (large -| O2);
					\legend{}
					\nextgroupplot[
					ybar, 
					bar width=22pt, 
					symbolic x coords={Strategy\#1,Strategy\#4},
					xtick=data, 
					xlabel={}, ylabel={Avg. runtime (\SI{}{\second})},
					xtick style={draw=none},
					xticklabels={Geom+IP, Pure IP},
					width=.5\textwidth,
					height=6cm,
					grid=major,
					ymin=0, 
					ymax=1000,
					ymode=log,
					log origin=infty,
					point meta=rawy,
					ytick={0,1,10,100},
					nodes near coords, 
					every node near coord/.append style={anchor=north, yshift=20pt}, 
					enlarge x limits=0.5, 
					legend style={at={(1.1,.82)}, anchor=west,legend columns=1}, 
					]
					%
					\addplot+[error bars/.cd, y dir=both, y explicit, error mark=triangle*] 
					table[row sep=\\,  y error plus=Max, y error minus=Min] {
						Category   Mean Min Max \\
						Strategy\#1  5.7600445    0.43536155  0.78675182 \\
						Strategy\#4  0.58286264   0.15821507  0.23533962 \\
					};
					\addlegendentry{Small} 
					\addplot+[error bars/.cd, y dir=both, y explicit, error mark=triangle*] 
					table[row sep=\\,  y error plus=Max, y error minus=Min] {
						Category   Mean Min Max \\
						Strategy\#1  16.07999322  3.10891244    8.73208025 \\
						Strategy\#4  19.46538837  6.22724802    7.76698391 \\
					};
					\addlegendentry{Medium} 
					\addplot+[error bars/.cd, y dir=both, y explicit, error mark=triangle*] 
					table[row sep=\\,  y error plus=Max, y error minus=Min] {
						Category   Mean Min Max \\
						Strategy\#1  51.65571241  5.39982014  7.67713518\\
						Strategy\#4  308.56465859  39.76806898  40.87631326\\
					};
					\addlegendentry{Large} 
				\end{groupplot}
			\end{tikzpicture}
		\caption{Performance of the different strategies on small, medium and large triangulations. The average, minimum and maximum of the found bicliques and runtimes are shown on the left and right hand side plot (resp.) for each strategy, the latter on a logarithmic scale.\label{fig:avgBiclique}}
		\end{figure}

		\subsection{Size of MILP formulations for higher dimensional simplicial partitions}\label{app:higher-dim-size}
		We constructed the known MILP formulations (DLog, Inc, MC, DCC, CC and GIB) for randomly chosen simplicial partitions in $\R^3$ and $\R^4$ of increasing size in terms of number of ground points and simplices.
		Table~\ref{tab:lpsize3d} shows a comparison of the number of rows, columns and binaries of the different formulations, as well as the number of simplices ($|\mathcal{P}|$), the number of ground points ($|V|$), the number of rank-2 edges $(\nu)$, and the number of edges of rank at least 3 $(\mu)$ of $\Hc_{\mathcal{S}}$, the number of edges ($\beta$)  of $\Hb_{\mathcal{S}}$, and the number of colors ($q$) and bicliques ($|\mathcal{B}|$) used in the GIB formulation.
		The instance names encode the dimension and the number of ground points of the simplicial partitions, and an ``r'' suffix indicates if the instance has been modified by the rank reduction procedure of Section~\ref{sec:rank_reduction}.
		Hence, every other row of the table shows the formulation size of the reduced version of the previous instance.
		For the instances that belong to rows with N/A, the construction of the blocking hypergraph failed due to running out of memory.
		The number of binary variables used by GIB is relatively close to those used by DLog (compared to the number of binaries in the othe formulations), while the number of continuous variables in DLog and the other formulations far exceeds the number of continuous variables of GIB.
		Also, observe that the rank reduction algorithm eliminated all conflicts of rank at least 3 in all instances.

		\section{Results on blocking triangulations for Section~\ref{sec:experiment:random}}\label{app:blocking}
		In this section we present computational results with  blocking triangulations for Section~\ref{sec:experiment:random} of the main article. 
		The two methods described in Section~\ref{sec:high-rank-conflict} were applied to resolve conflicts of rank 3.
		The sizes of the different triangulations and the corresponding formulations are summarized in Table~\ref{tab:randomSize} and Table~\ref{tab:lpsize}, respectively.
		The runtime of the models on the triangulations that were transformed by triangle subdivision as described in Section~\ref{sec:high-rank-conflict} are summarized in Table~\ref{tab:adjusted}, while Table~\ref{tab:unadjusted_incorrect} shows the results on blocking triangulation with the GIB formulation.
		A branch-and-cut method was implemented based on the network flow representation of the coloring constraints, based on~\citep{dobrovoczki_facet_2024}, see Remark~\ref{rem:network}.
		Column GIB+cuts shows the average runtime for the branch-and-cut method.
		A * sign marks the best average runtime for each scenario--triangulation size pair of Tables~\ref{tab:adjusted} and~\ref{tab:unadjusted_incorrect}.
		The triangle subdivision method was the best in most of the cases, due that the subdivision increased the number of triangles by at most 4.
		When the coloring constraints~(\ref{eq:colormip}) were used instead of the triangle subdivision, eight times out of nine only three colors were enough, and one time six colors were needed to properly color the blocking hypergraph.
		In this case, the branch-and-cut method outperformed by a little the GIB formulation on the medium and large instances, while on the small instances the model GIB dominated.
		
		\begin{table}
		\centering
		\small
		\caption{Runtime (\SI{}{\second}) of models on  random triangulations with triangle subdivisions.	}\label{tab:adjusted}
		\footnotesize
		\pgfplotstabletypeset[
		col sep=comma,
		header=true,
		columns/Month/.style={string type, column name=Month, column type=l},
		columns/Size/.style={string type, column name=Size, column type=r},
		columns/instances/.style={column name=\#instances, column type=r},
		columns/Dlog/.style={column name=DLog, fixed, precision=2, column type=r, fixed zerofill=true},
		columns/Inc/.style={column name=Inc, fixed, precision=2, column type=r, fixed zerofill=true},
		columns/MC/.style={column name=MC, fixed, precision=2, column type=r, fixed zerofill=true},
		columns/DCC/.style={column name=DCC, fixed, precision=2, column type=r, fixed zerofill=true},
		columns/CC/.style={column name=CC, fixed, precision=2, column type=r, fixed zerofill=true},
		columns/GIB/.style={column name=Subdiv.~+~IB, fixed, precision=2, column type=r, fixed zerofill=true},
		every head row/.style={before row={\toprule\hline}, after row=\midrule},
		every row no 3/.style={before row=\midrule},
		every row no 6/.style={before row=\midrule},
		every row no 9/.style={before row=\midrule},
		every last row/.style={after row=\bottomrule},
		every row 0 column 8/.style={postproc cell content/.style={@cell content/.add={$\bf}{$}}},
		every row 2 column 8/.style={postproc cell content/.style={@cell content/.add={$\bf}{^*$}}},
		every row 3 column 8/.style={postproc cell content/.style={@cell content/.add={$\bf}{$}}},
		every row 4 column 8/.style={postproc cell content/.style={@cell content/.add={$\bf}{$}}},
		every row 5 column 8/.style={postproc cell content/.style={@cell content/.add={$\bf}{^*$}}},
		every row 6 column 8/.style={postproc cell content/.style={@cell content/.add={$\bf}{^*$}}},
		every row 7 column 8/.style={postproc cell content/.style={@cell content/.add={$\bf}{^*$}}},
		every row 8 column 8/.style={postproc cell content/.style={@cell content/.add={$\bf}{$}}},
		every row 9 column 8/.style={postproc cell content/.style={@cell content/.add={$\bf}{^*$}}},
		opt/.list={3,4,5,6,7},
		boldoptstar/.list={8},
		row sep=\\,
		]{
			Month,Size,instances,Dlog,Inc,MC,DCC,CC,GIB\\
			April,small,2,3600,3600,3027.655,1810.73,1827.215,34.13\\
			,medium,1,3600,3600,3600,3600,3600,3600\\
			,large,6,3600,3600,3600,3600,3600,157.1683333\\
			June,small,2,3317.195,1855.905,2243.92,1811.525,2016.66,1806.315\\
			,medium,1,841.2,60.34,297.95,95.37,71.65,13.58\\
			,large,6,2094.458571,1965.85,1350.867143,1173.305714,1089.861429,14.765\\
			December,small,2,1875.53,3600,1828.86,1808.265,1816.735,1584.915\\
			,medium,1,3600,3600,3600,3600,3600,91.44\\
			,large,6,3600,3509.661667,3600,3251.416667,3184.771667,1491.113333\\
			Average,,,2981.850714,2957.40667,2655.636786,2409.989643,2390.416429,760.5240741\\
			Optimal,,27,6,6,10,10,10,23\\
		}
		
		\end{table}

		\begin{table}
		\centering
		\caption{Runtime (\SI{}{\second}) of models on  random triangulations with blocking triangles.\label{tab:unadjusted_incorrect}}
		\footnotesize
		\pgfplotstabletypeset[
		col sep=comma,
		header=true,
		columns/Month/.style={string type, column name=Month, column type=l},
		columns/Size/.style={string type, column name=Size, column type=r},
		columns/inst/.style={string type, column name=\#instances, column type=r},
		columns/DLog/.style={column name=DLog, fixed, precision=2, column type=r, fixed zerofill=true},
		columns/Inc/.style={column name=Inc, fixed, precision=2, column type=r, fixed zerofill=true},
		columns/MC/.style={column name=MC, fixed, precision=2, column type=r, fixed zerofill=true},
		columns/DCC/.style={column name=DCC, fixed, precision=2, column type=r, fixed zerofill=true},
		columns/CC/.style={column name=CC, fixed, precision=2, column type=r, fixed zerofill=true},
		columns/GIB/.style={column name=GIB, fixed, precision=2, column type=r, fixed zerofill=true},
		columns/GIB2/.style={column name=GIB+cuts, fixed, precision=2, column type=r, fixed zerofill=true},
		every head row/.style={before row={\toprule\hline}, after row=\midrule},
		every row no 3/.style={before row=\midrule},
		every row no 6/.style={before row=\midrule},
		every row no 9/.style={before row=\midrule},
		every last row/.style={after row=\bottomrule},
		every row 0 column 8/.style={postproc cell content/.style={@cell content/.add={$\bf}{^*$}}},
		every row 2 column 9/.style={postproc cell content/.style={@cell content/.add={$\bf}{$}}},
		every row 3 column 7/.style={postproc cell content/.style={@cell content/.add={$\bf}{^*$}}},
		every row 4 column 9/.style={postproc cell content/.style={@cell content/.add={$\bf}{^*$}}},
		every row 5 column 9/.style={postproc cell content/.style={@cell content/.add={$\bf}{$}}},
		every row 6 column 8/.style={postproc cell content/.style={@cell content/.add={$\bf}{$}}},
		every row 7 column 9/.style={postproc cell content/.style={@cell content/.add={$\bf}{$}}},
		every row 8 column 9/.style={postproc cell content/.style={@cell content/.add={$\bf}{^*$}}},
		every row 9 column 9/.style={postproc cell content/.style={@cell content/.add={$\bf}{$}}},
		opt/.list={3,4,5,6,7,8},
		boldopt/.list={9},
		row sep=\\,
		]{
			Month,Size,inst,DLog,Inc,MC,DCC,CC,GIB,GIB2\\
			April,small,2,1881.98,3223.025,2098.61,1809.545,1812.05,30.19,172.09\\
			,medium,1,3600,3600,3600,3600,3600,3600,3600\\
			,large,6,3600,3600,3178.226667,3355.258333,3600,1511.031667,1504.495\\
			June,small,2,2047.34,1866.69,1962.72,1813.49,1706.77,1805.9,1807.29\\
			,medium,1,260.48,63.29,866.26,151.69,179.27,8.16,6.89\\
			,large,6,2455.391667,2455.39,2450.195,1186.265,1523.938333,274.4016667,70.03833333\\
			December,small,2,1938.895,3600,1825.355,1809.425,1808.64,1804.49,1808.095\\
			,medium,1,3600,3600,3600,3600,3600,445.17,442.91\\
			,large,6,3600,3459.924,3600,3209.285,3561.891667,1181.386667,781.0183333\\
			Average,,,2856.821852,3027.207,2785.659259,2397.099259,2598.225556,1079.092222,954.0040741\\
			Optimal,,27,6,5,8,13,10,21,22\\
		}
		\end{table}
		
		\renewcommand{\thetable}{B.2}
		\begin{sidewaystable}[h]	
		\setlength{\tabcolsep}{3pt}	
		\centering
		\caption{Size of MILP formulation of the models on the random simplicial partitions in $\mathbb{R}^3$ and $\mathbb{R}^4$ of different sizes.}\label{tab:lpsize3d}
		\footnotesize
		\pgfplotstabletypeset[
		col sep=comma,
		header=true,
		empty cells with={N/A},
		columns/d/.style={column name=$d$, column type=r},
		columns/inst/.style={string type, column name=\#inst, column type=l},
		columns/T/.style={column name=$|\mathcal{P}|$, fixed, precision=0, column type=r, fixed zerofill=true},
		columns/V/.style={column name=$|V|$, fixed, precision=0, column type=r, fixed zerofill=true},
		columns/D/.style={column name=$\nu$, fixed, precision=0, column type=r, fixed zerofill=true},
		columns/H/.style={column name=$\mu$, fixed, precision=0, column type=r, fixed zerofill=true},
		columns/B/.style={column name=$\beta$, fixed, precision=0, column type=r, fixed zerofill=true},
		columns/ncolor/.style={column name=$q$, fixed, precision=0, column type=r, fixed zerofill=true},
		columns/nbiclique/.style={column name=$|\mathcal{B}|$, fixed, precision=0, column type=r, fixed zerofill=true},
		columns/gibrows/.style={column name=Rows, fixed, precision=0, column type=r, fixed zerofill=true},
		columns/gibcols/.style={column name=Cols, fixed, precision=0, column type=r, fixed zerofill=true},
		columns/gibbins/.style={column name=Bins, fixed, precision=0, column type=r, fixed zerofill=true},
		columns/dlogrows/.style={column name=Rows, fixed, precision=0, column type=r, fixed zerofill=true},
		columns/dlogcols/.style={column name=Cols, fixed, precision=0, column type=r, fixed zerofill=true},
		columns/dlogbins/.style={column name=Bins, fixed, precision=0, column type=r, fixed zerofill=true},
		columns/incrows/.style={column name=Rows, fixed, precision=0, column type=r, fixed zerofill=true},
		columns/inccols/.style={column name=Cols, fixed, precision=0, column type=r, fixed zerofill=true},
		columns/incbins/.style={column name=Bins, fixed, precision=0, column type=r, fixed zerofill=true},
		columns/mcrows/.style={column name=Rows, fixed, precision=0, column type=r, fixed zerofill=true},
		columns/mccols/.style={column name=Cols, fixed, precision=0, column type=r, fixed zerofill=true},
		columns/mcbins/.style={column name=Bins, fixed, precision=0, column type=r, fixed zerofill=true},
		columns/dccrows/.style={column name=Rows, fixed, precision=0, column type=r, fixed zerofill=true},
		columns/dcccols/.style={column name=Cols, fixed, precision=0, column type=r, fixed zerofill=true},
		columns/dccbins/.style={column name=Bins, fixed, precision=0, column type=r, fixed zerofill=true},
		columns/ccrows/.style={column name=Rows, fixed, precision=0, column type=r, fixed zerofill=true},
		columns/cccols/.style={column name=Cols, fixed, precision=0, column type=r, fixed zerofill=true},
		columns/ccbins/.style={column name=Bins, fixed, precision=0, column type=r, fixed zerofill=true},
		every head row/.style={
			before row={
				\toprule
				\hline
				&&&&&&&&
				& \multicolumn{3}{c@{\hspace{1em}}}{DLog}
				& \multicolumn{3}{c@{\hspace{1em}}}{Inc}
				& \multicolumn{3}{c@{\hspace{1em}}}{MC}
				& \multicolumn{3}{c@{\hspace{1em}}}{DCC}
				& \multicolumn{3}{c@{\hspace{1em}}}{CC}
				& \multicolumn{3}{c@{\hspace{1em}}}{GIB}\\
				\cmidrule(lr{1em}){10-12}
				\cmidrule(r{1em}){13-15}
				\cmidrule(r{1em}){16-18}
				\cmidrule(r{1em}){19-21}
				\cmidrule(r{1em}){22-24}
				\cmidrule(r{1em}){25-27}
			},
			after row=\midrule
		},
		every row no 2/.style={before row=\midrule},
		every row no 4/.style={before row=\midrule},
		every row no 6/.style={before row=\midrule},
		every row no 8/.style={before row=\midrule},
		every row no 10/.style={before row=\midrule},
		every row no 12/.style={before row=\midrule},
		every row no 14/.style={before row=\midrule},
		every row no 16/.style={before row=\midrule},
		every row no 18/.style={before row=\midrule},
		every row no 20/.style={before row=\midrule},
		every row no 22/.style={before row=\midrule},
		every row no 24/.style={before row=\midrule},
		every row no 26/.style={before row=\hline\midrule},
		every row no 28/.style={before row=\midrule},
		every row no 30/.style={before row=\midrule},
		every row no 32/.style={before row=\midrule},
		every last row/.style={after row=\bottomrule},
		row sep=\\,
		]{
			d,inst,T,V,D,H,B,ncolor,nbiclique,dlogrows,dlogcols,dlogbins,incrows,inccols,incbins,mcrows,mccols,mcbins,dccrows,dcccols,dccbins,ccrows,cccols,ccbins,gibrows,gibcols,gibbins\\
			3,{\ttfamily d3v11},22,11,17,8,352,6,6,15,98,10,45,87,21,101,88,22,27,110,22,17,33,22,27,23,12\\
			3,{\ttfamily d3v11r},30,14,40,0,0,0,7,15,130,10,61,119,29,136,120,30,35,150,30,20,44,30,32,21,7\\
			3,{\ttfamily d3v15},46,15,39,21,5515,15,9,17,196,12,93,183,45,201,184,46,51,230,46,21,61,46,37,39,24\\
			3,{\ttfamily d3v15r},80,25,187,0,0,0,12,19,334,14,161,319,79,347,320,80,85,400,80,31,105,80,53,37,12\\
			3,{\ttfamily d3v20},70,20,95,14,8219,11,10,19,294,14,141,279,69,302,280,70,75,350,70,26,90,70,44,41,21\\
			3,{\ttfamily d3v20r},103,29,267,0,0,0,13,19,426,14,207,411,102,443,412,103,108,515,103,35,132,103,59,42,13\\
			3,{\ttfamily d3v23},89,23,136,33,30196,21,12,19,370,14,179,355,88,381,356,89,94,445,89,29,112,89,51,56,33\\
			3,{\ttfamily d3v23r},141,37,480,0,0,0,15,21,580,16,283,563,140,603,564,141,146,705,141,43,178,141,71,52,15\\
			3,{\ttfamily d3v28},119,28,226,37,54024,17,14,19,490,14,239,475,118,506,476,119,124,595,119,34,147,119,60,59,31\\
			3,{\ttfamily d3v28r},188,46,793,0,0,0,18,21,768,16,377,751,187,800,752,188,193,940,188,52,234,188,86,64,18\\
			3,{\ttfamily d3v33},143,33,347,31,48141,14,14,21,588,16,287,571,142,607,572,143,148,715,143,39,176,143,65,61,28\\
			3,{\ttfamily d3v33r},216,52,1050,0,0,0,18,21,880,16,433,863,215,918,864,216,221,1080,216,58,268,216,92,70,18\\
			3,{\ttfamily d3v36},167,36,422,44,84775,14,15,21,684,16,335,667,166,706,668,167,172,835,167,42,203,167,70,65,29\\
			3,{\ttfamily d3v36r},262,61,1499,0,0,0,18,23,1066,18,525,1047,261,1111,1048,262,267,1310,262,67,323,262,101,79,18\\
			3,{\ttfamily d3v40},176,40,559,31,53102,9,16,21,720,16,353,703,175,746,704,176,181,880,176,46,216,176,76,65,25\\
			3,{\ttfamily d3v40r},248,58,1340,0,0,0,18,21,1008,16,497,991,247,1052,992,248,253,1240,248,64,306,248,98,76,18\\
			3,{\ttfamily d3v48},243,48,832,73,185269,19,18,21,988,16,487,971,242,1022,972,243,248,1215,243,54,291,243,88,85,37\\
			3,{\ttfamily d3v48r},415,89,3404,0,0,0,22,23,1678,18,831,1659,414,1751,1660,415,420,2075,415,95,504,415,137,111,22\\
			3,{\ttfamily d3v53},269,53,1051,70,215711,16,17,23,1094,18,539,1075,268,1131,1076,269,274,1345,269,59,322,269,91,86,33\\
			3,{\ttfamily d3v53r},434,93,3744,0,0,0,22,23,1754,18,869,1735,433,1831,1736,434,439,2170,434,99,527,434,141,115,22\\
			3,{\ttfamily d3v60},310,60,1395,84,262131,15,19,23,1258,18,621,1239,309,1302,1240,310,315,1550,310,66,370,310,102,94,34\\
			3,{\ttfamily d3v60r},489,105,4859,0,0,0,25,23,1974,18,979,1955,488,2063,1956,489,494,2445,489,111,594,489,159,130,25\\
			3,{\ttfamily d3v68},369,68,1836,100,374291,17,21,23,1494,18,739,1475,368,1546,1476,369,374,1845,369,74,437,369,114,106,38\\
			3,{\ttfamily d3v68r},617,128,7375,0,0,0,32,25,2488,20,1235,2467,616,2598,2468,617,622,3085,617,134,745,617,196,160,32\\
			3,{\ttfamily d3v75},400,75,2295,89,358598,15,20,23,1618,18,801,1599,399,1677,1600,400,405,2000,400,81,475,400,119,110,35\\
			3,{\ttfamily d3v75r},609,126,7132,0,0,0,30,25,2456,20,1219,2435,608,2564,2436,609,614,3045,609,132,735,609,190,156,30\\
			4,{\ttfamily d4v18},68,18,62,30,43860,25,6,20,354,14,137,339,67,360,340,68,74,408,68,25,86,68,35,49,31\\
			4,{\ttfamily d4v18r},304,46,721,0,0,0,18,24,1538,18,609,1519,303,1568,1520,304,310,1824,304,53,350,304,87,64,18\\
			4,{\ttfamily d4v20},102,20,76,41,197539,30,12,20,524,14,205,509,101,532,510,102,108,612,102,27,122,102,49,62,42\\
			4,{\ttfamily d4v20r},444,57,1173,0,0,0,20,24,2238,18,889,2219,443,2279,2220,444,450,2664,444,64,501,444,102,77,20\\
			4,{\ttfamily d4v22},142,22,91,55,,,,22,726,16,285,709,141,734,710,142,148,852,142,29,164,142,,,\\
			4,{\ttfamily d4v22r},671,74,2115,0,0,0,24,26,3375,20,1343,3354,670,3431,3355,671,677,4026,671,81,745,671,127,98,24\\
			4,{\ttfamily d4v24},172,24,115,57,,,,22,876,16,345,859,171,886,860,172,178,1032,172,31,196,172,,,\\
			4,{\ttfamily d4v24r},498,62,1432,0,0,0,23,24,2508,18,997,2489,497,2554,2490,498,504,2988,498,69,560,498,113,85,23\\
		}		
		\end{sidewaystable}

		\bibliographystyle{apalike}
		\bibliography{bibliography}
		\end{document}